\documentclass[11 pt,reqno]{amsart}
\usepackage{amssymb}
\usepackage{amsmath}

\input{diagrams}

\newcommand {\Ad}{\operatorname{Ad}}

\newcommand {\cor}[1]{\alpha_{#1}^{\vee}}	
\newcommand {\so}{\mathfrak{so}}
\newcommand {\wt}{\widetilde}

\newcommand {\IC}{\mathbb{C}}
\newcommand {\IN}{\mathbb{N}}                          
\newcommand {\IR}{\mathbb{R}}

\newcommand {\g}{\mathfrak{g}}
\newcommand {\A}{\mathcal A}
\newcommand {\D}{\mathcal D}
\newcommand {\V}{\mathcal V}

\newcommand {\half}[1]{\frac{#1}{2}}
\newcommand {\sqnm}[1]{\<#1,#1\>}

\newcommand {\End}{\operatorname{End}}
\newcommand {\Ker}{\operatorname{Ker}}
\newcommand {\<}{\langle}
\renewcommand {\>}{\rangle}

\newtheorem {theorem}{Theorem}[section]
\newtheorem {proposition}[theorem]{Proposition}
\newtheorem {corollary}[theorem]{Corollary}
\newtheorem {lemma}[theorem]{Lemma}
\numberwithin {equation}{section}
\renewcommand {\proof}{{\sc Proof.}\ }
\newcommand {\remark}{{\sc Remark.\ }}

\newcommand {\halmos}{$\Diamond$}

\newcommand {\KZ}{Knizhnik--Zamolodchikov  }
\newcommand {\ie}{{\it i.e., }}
\renewcommand {\halmos}{$\blacksquare$}

\newcommand {\Uh}{U_{\hbar}}
\newcommand {\Uhg}{\Uh\g}
\newcommand {\gl}{\mathfrak{gl}}
\renewcommand {\sl}{\mathfrak{sl}}
\newcommand {\Ugl}[1]{U\gl_{#1}}

\newcommand {\Uhgl}[1]{U_{\hbar}\gl_{#1}}
\newcommand {\Uhsl}[1]{U_{\hbar}\sl_{#1}}

\newcommand {\Smu}[1]{S^{\mu_{#1}}\IC^{k}}
\newcommand {\Shmu}[1]{S^{\mu_{#1}}_{\hbar}\IC^{k}}
\newcommand {\Sd}{\mathcal S^{d}\IC^{k}}
\newcommand {\Shd}{\mathcal S_{\hbar}^{d}\IC^{k}}
\newcommand {\fml}{[\negthinspace[\hbar]\negthinspace]}
\newcommand {\fmll}{[\negthinspace[h]\negthinspace]}
\newcommand {\ICh}{\IC\fml}
\newcommand {\young}[2]{\left(\substack{#1\\#2}\right)}

\newcommand {\M}{\mathcal M}

\newcommand {\SN}[1]{\mathfrak{S}_{#1}}
\renewcommand {\k}{^{(k)}}
\newcommand {\n}{^{(n)}}
\newcommand {\IY}{\mathbb Y}

\newcommand {\MM}[2]{\M_{#1,#2}}
\newcommand {\SM}{\mathcal S(\M^{*}_{k,n})}
\newcommand {\SMd}{\mathcal S^{d}(\M^{*}_{k,n})}
\newcommand {\SMh}{\mathcal S_{\hbar}(\M^{*}_{k,n})}
\newcommand {\SMhd}{\mathcal S_{\hbar}^{d}(\M^{*}_{k,n})}
\newcommand {\SMhdp}[2]{\mathcal S_{\hbar}^{d}(\M_{#1,#1})}

\newcommand {\h}{\mathfrak h}
\newcommand {\reg}{_{\operatorname{reg}}}
\newcommand {\hreg}{\h\reg}

\newcommand {\Bg}{B_{\g}}

\newcommand {\Ug}{U\g}

\renewcommand {\P}{\mathcal P}
\renewcommand {\dots}[1]{#1\cdots#1}

\renewcommand {\SS}{\mathfrak S}
\newcommand {\tr}{\operatorname{tr}}

\newcommand {\bin}[2]
{\begin{bmatrix}#1\\#2\end{bmatrix}}
\newcommand {\eg}{{\it e.g., }}
\newcommand {\Sh}[2]{\mathcal S_{\hbar}(\mathcal M_{#1,#2}^{*})}
\newcommand {\XX}{\mathbf X}
\newcommand {\res}{^{\circ}}

\newcommand {\p}{^{(p)}}
\newcommand {\SmC}[1]{\mathcal S^{#1}\IC^{k}}

\newcommand {\hbbar}{\overline{h}}
\newcommand {\KKZ}{_{\scriptscriptstyle{\operatorname{KZ}}}}
\newcommand {\back}{\!\!\!\!}
\newcommand {\e}{\mathfrak e}
\renewcommand {\sp}{\mathfrak{sp}}
\newcommand {\id}{\operatorname{id}}
\newcommand {\SmCh}[1]{\mathcal S_{\hbar}^{#1}\IC^{k}}
\newcommand {\nablak}{\nabla_{\kappa}}
\newcommand {\lldots}{\ldots\negthinspace}

\begin{document}

\title
[a Kohno--Drinfeld theorem for Quantum Weyl Groups]
{a Kohno--Drinfeld theorem for Quantum Weyl Groups}
\author
[V. Toledano Laredo]
{Valerio Toledano Laredo}
\address{
MSRI                                   \newline
1000 Centennial Drive                  \newline
Berkeley, CA 94720--5070               \newline
toledano@msri.org		       \newline
                                       \newline
permanent address :                    \newline
                                       \newline
Institut de Mathematiques de Jussieu   \newline
UMR 7586, Case 191                     \newline
Universite Pierre et Marie Curie       \newline
4, Place Jussieu                       \newline
F--75252 Paris Cedex 05                \newline
toledano@math.jussieu.fr}
\date{8 May 2001}
\begin{abstract}
Let $\g$ be a complex, simple Lie algebra with Cartan
subalgebra $\h$ and Weyl group $W$. In \cite {MTL},
we introduced a new, $W$--equivariant flat connection
on $\h$ with simple poles along the root hyperplanes
and values in any finite--dimensional $\g$--module $V$.
It was conjectured in \cite{TL} that its monodromy is
equivalent to the quantum Weyl group action of the
generalised braid group of type $\g$ on $V$ obtained
by regarding the latter as a module over the quantum
group $\Uhg$. In this paper, we prove this conjecture
for $\g=\sl_{n}$.
\end{abstract}
\maketitle

\section{Introduction}
%---------------------

One of the many virtues of quantum groups is their ability
to describe the monodromy of certain first order systems
of Fuchsian PDEs. If $\Uhg$ is the Drinfeld--Jimbo quantum
group of the complex, simple Lie algebra $\g$, the
universal $R$--matrix of $\Uhg$ yields a representation
of Artin's braid group on $n$ strings $B_{n}$ on the $n$--fold
tensor product $V^{\otimes n}$ of any finite--dimensional
$\Uhg$--module $V$. A fundamental, and paradigmatic result
of Kohno and Drinfeld establishes the equivalence of this
representation with the monodromy of the \KZ equations for
$\g$ with values in $V^{\otimes n}$ \cite{Dr3,Dr4,Dr5,Ko1}.
Lusztig, and independently Kirillov--Reshetikhin and
Soibelman realised that $\Uhg$ also yields representations
of another braid group, namely the generalised braid
group $\Bg$ of Lie type $\g$ \cite{Lu1,KR,So}. Whereas
the $R$--matrix representation is a deformation of the
natural action of the symmetric group ${\mathfrak S}_{n}$
on $n$--fold tensor products, these representations of
$\Bg$ deform the action of (a finite extension of) the
Weyl group $W$ of $\g$ on any finite--dimensional $\g$--module
V.\\

The aim of this paper is to show that these quantum Weyl
group representations describe the monodromy of the flat
connection introduced in \cite{MTL} and, independently,
in \cite{FMTV}. More precisely, realise $\Bg$ as the
fundamental group of the orbit space $\hreg/W$ of the set of
regular elements of a Cartan subalgebra $\h$ of $\g$ under
the action of $W$ \cite{Br}. Then, one can define a flat
vector bundle $(\V,\nabla_{\kappa})$ with fibre $V$ over
$\hreg/W$ \cite{MTL}. The connection $\nabla_{\kappa}$
depends upon a parameter $\hbar\in\IC$ and it was conjectured
in \cite{TL} that, when $\hbar$ is regarded as a formal
variable, its monodromy is equivalent to the quantum Weyl
group action of $\Bg$ on $V$. This conjecture was checked
in \cite{TL} for a number of pairs $(\g,V)$ including vector
representations of classical Lie algebras and adjoint
representations of simple Lie algebras.\\

In the present paper, we prove this conjecture for $\g=\sl_{n}$,
so that $\Bg=B_{n}$. The proof relies on the Kohno--Drinfeld
theorem for $\Uhsl{k}$ via the use of the dual pair
$(\gl_{k},\gl_{n})$. Our main observation is that the
duality between $\gl_{k}$ and $\gl_{n}$ derived from their
joint action on the space $\M_{k,n}$ of $k\times n$ matrices
exchanges $\nabla_{\kappa}$ for $\sl_{n}$ and the \KZ connection
for $\sl_{k}$, thus acting as a simple--minded integral
transform. This shows the equivalence of the monodromy
representation of $\nabla_{\kappa}$ for $\sl_{n}$ with
a suitable $R$--matrix representation for $\Uh\sl_{k}$.
The proof is completed by noting that the duality
between $\Uh\gl_{k}$ and $\Uh\gl_{n}$ exchanges the
$R$--matrix representation of $\Uhsl{k}$  with the
quantum Weyl group representation of $\Uhsl{n}$.

This may be schematically summarised by the following
diagram

\newarrow{Corresponds}{<}{-}{-}{-}{>}
\begin{equation*}
\begin{diagram}[height=4em,width=6em]
\nabla\KKZ,\sl_{k}&\rCorresponds^{\M_{k,n}}  &\nabla_{\kappa},\sl_{n}\\
\dCorresponds^{KD}&                          &\dCorresponds\\
R^{\vee},\Uhsl{k} &\rCorresponds_{\M^{\hbar}_{k,n}}&W_{\hbar},\Uhsl{n}
\end{diagram}
\end{equation*}

The structure of the paper is as follows. In section
\ref{se:flat}, we give the construction of the connection
$\nabla_{\kappa}$ following \cite{MTL}. We show in
section \ref{se:KZ Howe} that the duality between
$\gl_{k}$ and $\gl_{n}$ identifies the \KZ connection
for $n$--fold tensor products of symmetric powers
of the vector representation of $\sl_{k}$ and the
connection $\nabla_{\kappa}$ for $\sl_{n}$. In
section \ref{se:Uhg} we recall the definition
of the Drinfeld--Jimbo quantum groups $\Uhgl{k}$
and $\Uhgl{n}$ and, in section \ref{se:q Howe} show
how they jointly act on the quantum $k\times n$
matrix space $\SMh$. The corresponding $R$--matrix
and quantum Weyl group representations of $B_{n}$
on $\SMh$ are shown to coincide in section \ref
{se:Bn on Mkn}. Section \ref{se:main} contains
our main result.\\

{\bf Acknowledgements. } This paper was begun at the
Reseach Institute for Mathematical Sciences of Kyoto
University. I am very grateful to M. Kashiwara for
his invitation to spend the summer of 1999 at RIMS
and to RIMS for its hospitality and financial support.
During my stay, I greatly benefitted from very stimulating
and informative discussions with M. Kashiwara and B.
Feigin. I also wish to express my gratitude to
A. D'Agnolo, P. Baumann, B. Enriquez, J. Millson,
R. Rouquier and P. Schapira for innumerable, useful
and friendly conversations.

\section{Flat connections on $\h\reg$}\label{se:flat}
%=====================================

The results in this section are due to J. Millson and
the author \cite{MTL}. They were obtained independently
by De Concini around 1995 (unpublished). Let $\g$ be
a complex, simple Lie algebra with Cartan subalgebra
$\h$ and root system $R\subset\h^{*}$. Let $\hreg=\h
\setminus\bigcup_{\alpha\in R}\Ker(\alpha)$ be the set
of regular elements in $\h$ and $V$ a finite--dimensional
$\g$--module. We shall presently define a flat connection
on the topologically trivial vector bundle $\hreg\times
V$ over $\hreg$. We need for this purpose the following
simple flatness criterion due to Kohno \cite{Ko2}. Let
$B$ be a complex, finite--dimensional vector space and
$\A=\{H_{i}\}_{i\in I}$ a finite collection of hyperplanes
in $B$ determined by the linear forms $\phi_{i}\in B^{*}$,
$i\in I$.

\begin{lemma}\label{le:kohno}
Let $F$ be a finite--dimensional vector space and
$\{r_{i}\}\subset\End(F)$ a family indexed by
$I$. Then,
\begin{equation}
\nabla=d-\sum_{i\in I}\frac{d\phi_{i}}{\phi_{i}}r_{i}
\end{equation}
defines a flat connection on $(B\setminus\A)\times F$
iff, for any subset $J\subseteq I$ maximal for the
property that $\bigcap_{j\in J}H_{j}$ is of codimension
2, the following relations hold for any $j\in J$
\begin{equation}
[r_{j},\sum_{j'\in J}r_{j'}]=0
\end{equation}
\end{lemma}

For any $\alpha\in R$, choose root vectors $e_{\alpha}
\in\g_{\alpha},f_{\alpha}\in\g_{-\alpha}$ such that
$[e_{\alpha},f_{\alpha}]=h_{\alpha}=\alpha^{\vee}$
and let
\begin{equation}
\kappa_{\alpha}=
\half{\sqnm{\alpha}}
(e_{\alpha}f_{\alpha}+f_{\alpha}e_{\alpha})
\in U\g
\end{equation}
be the truncated Casimir operator of the $\sl_{2}
(\IC)$--subalgebra of $\g$ spanned by $e_{\alpha},
h_{\alpha},f_{\alpha}$. Note that $\kappa_{\alpha}$
does not depend upon the particular choice of $e_{\alpha}$
and $f_{\alpha}$ and that $\kappa_{-\alpha}=\kappa
_{\alpha}$.

\begin{theorem}\label{th:casimir flat}
The one--form
\begin{equation}\label{eq:Casimir connection}
\nabla_{\kappa}^{h}=
d-h\sum_{\alpha\succ 0}
\frac{d\alpha}{\alpha}\kappa_{\alpha}=
d-\frac{h}{2}\sum_{\alpha\in R}
\frac{d\alpha}{\alpha}\kappa_{\alpha}
\end{equation}
defines, for any $h\in\IC$, a flat connection on
$\hreg\times V$.
\end{theorem}
\proof By lemma \ref{le:kohno}, we must prove
that for any rank 2 subsystem $R_{0}\subseteq R$,
the following holds for any $\alpha\in R_{0}^{+}=
R_{0}\cap R^{+}$
\begin{equation}\label{eq:rank 2}
[\kappa_{\alpha},\sum_{\beta\in R_{0}^{+}}\kappa_{\beta}]=0
\end{equation}
This may be proved by an explicit computation by
considering in turn the cases where $R_{0}$ is of
type $A_{1}\times A_{1}$, $A_{2}$, $B_{2}$ or
$G_{2}$ but is more easily settled by the following
elegant observation of A. Knutson \cite{Kn}. Let
$\g_{0}\subseteq\g$ be the semi--simple Lie algebra
with root system $R_{0}$, $\h_{0}\subset\h$ its Cartan
subalgebra and $C_{0}\in Z(\Ug_{0})$ its Casimir
operator. Then, $\sum_{\beta\in R_{0}^{+}}\kappa_{\beta}
-C_{0}$ lies in $U\h_{0}$ so that \eqref{eq:rank 2}
holds since $\kappa_{\alpha}$ commutes with $\h_{0}$
\halmos\\

Let $G$ be the complex, connected and simply--connected
Lie group with Lie algebra $\g$, $T$ its torus with Lie
algebra $\h$, $N(T)\subset G$ the normaliser of $T$ and
$W=N(T)/T$ the Weyl group of $G$. Let $\Bg=\pi_{1}(\hreg
/W)$ be the generalised braid group of type $\g$ and
$\sigma:\Bg\rightarrow N(T)$ a homomorphism compatible with
\begin{equation}\label{eq:extension}
\begin{diagram}[height=2.5em]
\Bg&\rTo^{\sigma}&N(T) \\
   &\rdTo        &\dTo \\
   &             &W
\end{diagram}
\end{equation}
\cite{Ti}. We regard $\Bg$ as acting on $V$ via $\sigma$.
Let $\wt{\hreg}\xrightarrow{p}\hreg$ be the universal
cover of $\hreg$ and $\hreg/W$.

\begin{proposition}\label{th:existence}
The one--form $p^{*}\nabla_{\kappa}^{h}$ defines a $\Bg$--equivariant
flat connection on $\wt\hreg\times V=p^{*}(\hreg\times V)$.
It therefore descends to a flat connection on the vector
bundle
\begin{equation}
\begin{diagram}[height=2em,width=2em]
V&\rTo&\wt{\hreg}\times_{\Bg}V \\
 &    &\dTo                    \\
 &    &\hreg/W
\end{diagram}
\end{equation}
which is reducible with respect to the weight space
decomposition of $V$ and unitary if $h\in i\IR$.
\end{proposition}
\proof
The action of $\Bg$ on $\Omega^{\bullet}(\wt\hreg,V)
=\Omega^{\bullet}(\wt\hreg)\otimes V$ is given by
$\gamma\rightarrow(\gamma^{-1})^{*}\otimes\sigma
(\gamma)$. Thus, if $\gamma\in\Bg$ projects onto
$w\in W$, we get using $p\cdot\gamma^{-1}=w^{-1}
\cdot p$,
\begin{equation}\label{eq:equivariant}
\gamma\medspace p^{*}\nabla_{\kappa}^{h}\medspace\gamma^{-1}=
d-
\half{h}\sum_{\alpha\in R}dp^{*}w\alpha/p^{*}w\alpha
\otimes\sigma(\gamma)\kappa_{\alpha}\sigma(\gamma)^{-1}
\end{equation}
Since $\kappa_{\alpha}=\half{\sqnm{\alpha}}(e_{\alpha}
f_{\alpha}+f_{\alpha}e_{\alpha})$ is independent
of the choice of the root vectors $e_{\alpha},
f_{\alpha}$, $\Ad(\sigma(\gamma))\kappa_{\alpha}
=\kappa_{w\alpha}$ and \eqref{eq:equivariant} is
equal to $p^{*}\nabla_{\kappa}^{h}$ as claimed.
$p^{*}\nabla_{\kappa}^{h}$ is flat by theorem \ref{th:casimir flat},
commutes with the fibrewise action of $\h$ because
each $\kappa_{\alpha}$ is of weight $0$ and is unitary
because the $\kappa_{\alpha}$ are self--adjoint \halmos\\

Thus, for any homomorphism $\sigma:\Bg\rightarrow N(T)$
compatible with \eqref{eq:extension}, proposition
\ref{th:existence} yields a monodromy representation
$\rho^{\sigma}_{h}:\Bg\rightarrow GL(V)$ which permutes the
weight spaces compatibly with $W$. By standard ODE
theory, $\rho^{\sigma}_{h}$ depends analytically on the complex
parameter $h$ and, when $h=0$, is equal to the action
of $\Bg$ on $V$ given by $\sigma$. We record for later
use the following elementary

\begin{proposition}\label{pr:recipe}
Let $\gamma\in\Bg=\pi_{1}(\hreg/W)$ and $\wt
\gamma:[0,1]\rightarrow\hreg$ be a lift of
$\gamma$. Then,
\begin{equation}
\rho^{\sigma}_{h}(\gamma)=\sigma(\gamma)\P(\wt\gamma)
\end{equation}
where $\P(\wt\gamma)\in GL(V)$ is the parallel
transport along $\wt\gamma$ for the connection
$\nabla_{\kappa}^{h}$ on $\hreg\times V$.
\end{proposition}
\proof
Let $\wt{\wt\gamma}:[0,1]\rightarrow\wt{\hreg}$
be a lift of $\gamma$ and $\wt{\gamma}$ so that
$\wt{\wt\gamma}(1)=\gamma^{-1}\wt{\wt\gamma}(0)$.
Then, since the connection on $p^{*}(\hreg\times
V)$ is the pull--back of $\nabla_{\kappa}^{h}$, and that
on $\left(p^{*}(\hreg\times V)\right)/\Bg$ the
quotient of $p^{*}\nabla_{\kappa}^{h}$, we find
\begin{equation}
\rho^{\sigma}_{h}(\gamma)=
\P(\gamma)=
\sigma(\gamma)\P(\wt{\wt\gamma})=
\sigma(\gamma)\P(\wt\gamma)
\end{equation}
\halmos\\

By \cite{Br}, $\Bg$ is presented on generators
$T_{i}$, $i=1\ldots n$ labelled by a choice of
simple roots $\alpha_{i}$ of $R$ with relations
\begin{equation}\label{eq:braid reln}
T_{i}T_{j}T_{i}\cdots=
T_{j}T_{i}T_{j}\cdots
\end{equation}
for any $1\leq i<j\leq n$ where each side of
\eqref{eq:braid reln} has a number of factors
equal to the order of $s_{i}s_{j}$ in $W$ and
$s_{k}\in W$ is the orthogonal reflection
across the hyperplane $\Ker(\alpha_{k})$.
$T_{i}$ projects onto $s_{i}\in W$. An
explicit choice of representatives of $T_{1},
\ldots,T_{n}$ in $\pi_{1}(\hreg/W)$ may be
given as follows.
Let $t\in\hreg$ lie in the fundamental Weyl
chamber so that $\<t,\alpha\>>0$ for any
$\alpha\in R_{+}$. Note that for any simple root
$\alpha_{i}$, the intersection $t_{\alpha_{i}}
=t-\half{1}\<t,\alpha_{i}\>\cor{i}$
of the affine line $t+\IC\cdot\alpha_{i}^
{\vee}$ with $\Ker(\alpha_{i})$ does not lie
in any other root hyperplane $\Ker(\beta)$,
$\beta\in R\setminus \{\alpha_{i}\}$. Indeed,
if $\<t_{\alpha_{i}},\beta\>=0$ then
\begin{equation}
2\<t,\beta\>=
\<t,\alpha_{i}\>\<\cor{i},\beta\>=
\<t,\beta-s_{i}\beta\>
\end{equation}
whence $\<t,\beta\>=-\<t,s_{i}\beta\>$,
a contradiction since $s_{i}$ permutes positive
roots different from $\alpha_{i}$. Let now $D$
be an open disc in $t+\IC\cdot\cor{i}$ of center
$t_{\alpha_{i}}$ such that its closure $\overline{D}$
does not intersect any root hyperplane other than
$\Ker(\alpha_{i})$. Consider the path $\gamma_{i}:
[0,1]\rightarrow t+\IC\cdot\cor{i}$ from $t$ to $s_{i}t$
determined by $\left.\gamma_{i}\right|_{[0,1/3]\cup[2/3,1]}$
is affine and lies in $t+\IR\cdot\alpha^{\vee}\setminus
D$, $\gamma_{i}(1/3),\gamma_{i}(2/3)\in\partial
\overline{D}$ and $\left.\gamma_{i}\right|_{[1/3,2/3]}$
is a semicircular arc in $\partial\overline{D}$,
positively oriented with respect to the natural
orientation of $t+\IC\cdot\alpha_{i}$. Then, the image
of $\gamma_{i}$ in $\hreg/W$ is a representative of
$T_{i}$ in $\pi_{1}(\hreg/W,Wt)$ \cite{Br}.

\section{Knizhnik--Zamolodchikov equations and dual pairs}\label{se:KZ Howe}
%=========================================================

We show in this section that the joint action of $\gl_{k}$
and $\gl_{n}$ on the space $\M_{k,n}$ of $k\times n$ matrices
identifies the connection $\nablak$ for $\g=\sl_{n}$ and the
\KZ connection for $\sl_{k}$. Let $\SM=\IC[x_{11},\ldots,x_{kn}]$
be the algebra of polynomial functions on $\M_{k,n}$. The group
$GL_{k}\times GL_{n}$ acts on $\SM$ by
\begin{equation}
(g_{k},g_{n})\medspace p(x)=
p(g_{k}^{t}xg_{n})
\end{equation}
and leaves the homogeneous components $\SMd$ of $\SM$,
$d\in\IN$, invariant.
The decomposition of $\SM$ under $GL_{k}\times GL_{n}$
is well--known (see \eg \cite[\S 132]{Zh} which we
follow closely or \cite[\S 1.4]{Mc} \footnote{I
am grateful to M. Vergne for pointing out that
the decomposition \eqref{eq:Howe} was known long
before \cite{Zh} and the work of Howe \cite{Ho} and to
M. Brion for providing the reference \cite{Mc}.}).
Let $N_{k},N_{n}$ be the groups of $k\times k$ and
$n\times n$ upper triangular unipotent matrices
respectively.

\begin{lemma}\label{le:inva}
\begin{equation}
\SM^{N_{k}\times N_{n}}=
\IC[\Delta_{1},\ldots,\Delta_{\min(k,n)}]
\end{equation}
where $\Delta_{l}(x)=\det(x_{ij})_{1\leq i,j\leq l}$
is the $l$th principal minor of the matrix $x$.
\end{lemma}
\proof
Assume for simplicity that $k\leq n$. Let $\D\subset
\SM$ be the subset of matrices $x$ such that $\Delta
_{i}(x)\neq 0$ for $i=1\ldots k$. By the Gauss
decomposition, any $x\in\D$ is conjugate under $N_{k}
^{t}\times N_{n}$ to a unique $k\times n$ matrix $d
(x)$ with the same principal minors as $x$, diagonal
principal $k\times k$ block and the remaining columns
equal to zero. Consider now the $k\times n$ matrix
\begin{equation}
m(x)=\left(\begin{array}{ccccccccc}
\Delta_{1}(x)&\Delta_{2}(x)&\Delta_{3}(x)&\cdots&\Delta_{k-1}(x)&\Delta_{k}(x)&0&\cdots&0 \\
           -1&            0&            0&\cdots&              0&            0&0&\cdots&0 \\
            0&           -1&            0&\cdots&              0&            0&0&\cdots&0 \\
%            0&            0&           -1&\cdots&              0&            0&0&\cdots&0 \\
       \vdots&       \vdots&       \vdots&\cdots&         \vdots&       \vdots&\vdots&\cdots&\vdots \\
            0&            0&            0&\cdots&             -1&            0&0&\cdots&0 \\
\end{array}\right)
\end{equation}
Since $\Delta_{i}(m(x))=\Delta_{i}(x)$, $1\leq i
\leq k$, $m(x)$ is also conjugate to $d(x)$, and
therefore to $x$, under $N_{k}^{t}\times N_{n}$.
Thus, by density of $\D$, a polyonomial $p\in\SM$
is invariant under $N_{k}\times N_{n}$ iff it is
a function of $m(x)$ and therefore iff it is a
polynomial in $\Delta_{1},\ldots,\Delta_{k}$
\halmos\\

Let $\IY_{p}\subset\IN^{p}$ be the set of Young
diagrams with at most $p$ rows. For $\lambda\in
\IY_{p}$, set $|\lambda|=\sum_{i=1}^{p}\lambda_{i}$
and let $V_{\lambda}^{(p)}$ be the irreducible
representation of $GL_{p}(\IC)$ of highest weight
$\lambda$.

\begin{theorem}\label{th:Howe}
As $GL_{k}\times GL_{n}$--modules,
\begin{equation}\label{eq:Howe}
\SMd\cong
\bigoplus_{\substack{\lambda\in\IY_{\min(k,n)},\\|\lambda|=d}}
V_{\lambda}^{(k)}\otimes V_{\lambda}^{(n)}
\end{equation}
\end{theorem}
\proof
Assume again $k\leq n$ for simplicity. By lemma
\ref{le:inva}, the highest weight vectors for
the action of $GL_{k}(\IC)\times GL_{n}(\IC)$
on $\SM$ are the polynomials in $\Delta_{1},\ldots,
\Delta_{k}$ which are eigenvectors for the torus
of $GL_{k}\times GL_{n}$. Since $\Delta_{l}$
is of weight $\varpi_{l}\k\oplus\varpi_{l}\n$,
where $\varpi_{l}^{(p)}$ %=\theta_{1}+\cdots+\theta_{l}
is the $l$th fundamental weight of $GL_{p}$, the
highest weight vectors are the monomials $\Delta_{1}
^{m_{1}}\cdots\Delta_{k}^{m_{k}}$ with corresponding
pair of Young diagrams $(\lambda,\lambda)$ where
\begin{equation}
\lambda=(m_{1}+\cdots+m_{k},m_{2}+\cdots+m_{k},\ldots,m_{k})
\end{equation}
Thus, \eqref{eq:Howe} holds since $\Delta_{l}$
is a homogeneous function of degree $l$ \halmos\\

As a $\gl_{k}$--module,
\begin{equation}
\SM=
\IC[x_{11},\ldots,x_{k1}]
\dots{\otimes}
\IC[x_{1n},\ldots,x_{kn}]
\end{equation}
and is therefore acted upon by the $\gl_{k}$--intertwiners
$\wt{\Omega}_{ij}\k$, $1\leq i<j\leq n$, defined by
\begin{equation}\label{eq:Omegaij}
\wt{\Omega}_{ij}\k=
\sum_{a}
1^{\otimes(i-1)}\otimes X_{a}\otimes 1^{\otimes(j-i-1)}
\otimes X^{a}\otimes 1^{\otimes(n-j)}
\end{equation}
where $\{X_{a}\}$, $\{X^{a}\}$ are dual basis of $
\gl_{k}$ with respect to the pairing $\<X,Y\>=\tr
(XY)$. On the other hand, as a $\gl_{n}$--module,
$\SM$ is acted upon by the operators $\kappa_{ij}
\n$, $1\leq i<j\leq n$, where
\begin{equation}\label{eq:Cij}
\kappa_{ij}\n=
e_{\alpha}f_{\alpha}+
f_{\alpha}e_{\alpha}
\end{equation}
is the truncated Casimir operator of the $\sl_{2}
(\IC)$--subalgebra of $\gl_{n}$ corresponding to
the root $\alpha=\theta_{i}-\theta_{j}$.
Let $e_{1},\ldots,e_{p}$ be the canonical basis of
$\IC^{p}$ and $E_{ab}^{(p)}e_{c}=\delta_{bc}e_{a}$,
$1\leq a,b\leq p$ the corresponding basis of $\gl_
{p}$ with dual basis $E_{ba}^{(p)}$. Let $1\leq i
<j\leq n$, then

\begin{proposition}\label{pr:Omegaij=Cij}
The following holds on $\SM$
\begin{equation}
2\wt{\Omega}_{ij}\k=\kappa_{ij}\n-E_{ii}\n-E_{jj}\n
\end{equation}
\end{proposition}
\proof By \eqref{eq:Omegaij}, $\wt{\Omega}_{ij}\k$
acts on $\SM$ as
\begin{equation}
\wt{\Omega}_{ij}\k=
\sum_{1\leq a,b\leq k}
x_{ai}\partial_{bi}x_{bj}\partial_{aj}
\end{equation}
where $x_{rc}$ and $\partial_{rc}$ are the operators of
multiplication by and derivation with respect to $x_{rc}$.
On the other hand, given that the $\sl_{2}(\IC)$--triple
$\{e_{\alpha},h_{\alpha},f_{\alpha}\}$ corresponding to
the root $\alpha=\theta_{i}-\theta_{j}$ of $\gl_{n}$ is
$\{E_{ij}\n,E_{ii}\n-E_{jj}\n,E_{ji}\n\}$, the following
holds on $\SM$
\begin{equation}
\kappa_{ij}\n=
\sum_{1\leq a,b\leq k}
x_{ai}\partial_{aj}x_{bj}\partial_{bi}+
x_{bj}\partial_{bi}x_{ai}\partial_{aj}
\end{equation}
Substracting, we find
\begin{equation}
2\wt{\Omega}_{ij}\k-\kappa_{ij}\n=-
\sum_{1\leq a,b\leq k}
\delta_{ab}x_{ai}\partial_{bi}+
\delta_{ab}x_{bj}\partial_{aj}=
-E_{ii}\n-E_{jj}\n
\end{equation}
as claimed \halmos\\

Let $\lambda\in\IY_{\min(k,n)}$ and $V_{\lambda}
\n$ the corresponding simple $GL_{n}$--module.
By theorem \ref{th:Howe}, $V_{\lambda}\n$ may
be identified with the subspace of vectors of
highest weight $\lambda$ for the action of
$\gl_{k}$ on $\SM$. Denote by $\iota:V_{\lambda}
\rightarrow\SM$ the corresponding $\gl_{n}
$--equivariant embedding and let $\mu=(\mu_{1},
\ldots,\mu_{n})\in\IN^{n}$ be a weight of
$V_{\lambda}\n$.

\begin{lemma}\label{le:singular=weight}
The embedding $\iota$ maps the subspace $V_{\lambda}
\n[\mu]\subset V_{\lambda}\n$ of weight $\mu$ onto
the subspace $M_{\lambda}^{\mu}$ of vectors of
highest weight $\lambda$ for the action of
$\gl_{k}$ on
\begin{equation}
\Smu{}=\Smu{1}\dots{\otimes}\Smu{n}
\subset
\IC[x_{11},\ldots,x_{k1}]\dots{\otimes}\IC[x_{1n},\ldots,x_{kn}]
\end{equation}
where $\Smu{j}$ is the space of polynomials in
$x_{1j},\ldots,x_{kj}$ which are homogeneous of
degree $\mu_{j}$. The corresponding isomorphism
\begin{equation}
\bigoplus_{\nu\in\SN{n}\mu}V_{\lambda}\n[\nu]
\cong
\bigoplus_{\nu\in\SN{n}\mu}M_{\lambda}^{\nu}
\end{equation}
is equivariant with respect to $\SN{n}$ which
acts on $\bigoplus_{\nu\in\SN{n}\mu}S^{\nu}\IC^{k}$
by permuting the tensor factors and on $V_{\lambda}$
by regarding $\SS_{n}$ as the subgroup of permutation
matrices of $GL_{n}(\IC)$.
\end{lemma}
\proof
The equality $\iota(V_{\lambda}\n[\mu])=M_{\lambda}
^{\mu}$ holds because $\Smu{}$ is the subspace of
$\SM$ of weight $\mu$ for the $\gl_{n}$--action
since $E\n_{ii}x_{rj}^{m}=\delta_{ij}m x_{rj}^{m}$.
The $\SS_{n}$--equivariance stems from the fact
that the permutation of the tensor factors in
$S^{\bullet}\IC^{k}\dots{\otimes}S^{\bullet}\IC^{k}\cong
S^{\bullet}(\IC^{k}\otimes\IC^{n})$ is given by the
action of $\SN{n}\subset GL_{n}(\IC)$ action on
$\IC^{n}$ \halmos\\

Let $\D_{n}=\{(z_{1},\ldots,z_{n})\in\IC^{n}|\medspace
z_{i}=z_{j}\medspace\medspace\text{for some $1\leq i<j
\leq n$}\}$ and $X_{n}=\IC^{n}\setminus\D_{n}$. Regard
$\IC^{n}_{0}=\{(z_{1},\ldots,z_{n})\in\IC^{n}|\medspace
\sum_{j=1}^{n}z_{j}=0\}$ as the Cartan subalgebra of
diagonal matrices in $\sl_{n}$ and $X_{n}^{0}=
\IC^{n}_{0}\setminus\D_{n}$ as the set of its regular
elements. Since the inclusion $X_{n}^{0}\subset X_{n}$
is a homotopy equivalence, $\pi_{1}(X_{n})\cong\pi_{1}
(X_{n}^{0})=B_{n}$ are generated by $T_{1}\ldots T_{n-1}$
with
\begin{alignat}{2}
T_{i}T_{j}&=T_{j}T_{i}\qquad& &\text{if $|i-j|\geq 2$}\\
T_{i}T_{i+1}T_{i}&=T_{i+1}T_{i}T_{i+1}\qquad&
&\text{$i=1\ldots n-1$}
\end{alignat}
Define $\Omega\k_{ij}\in\End_{\gl_{k}}(\Smu{})$ by
\eqref{eq:Omegaij} where now $\{X_{a}\},\{X^{a}\}$
are dual basis of $\sl_{k}$ and extend the connection
\eqref{eq:Casimir connection} to $X_{n}$ in the obvious
way. The following is the main result of this section.

\begin{theorem}\label{th:KZ=KZ}
$f:X_{n}\rightarrow M^{\mu}_{\lambda}\subset
\Smu{1}\dots{\otimes}\Smu{n}$ is a horizontal
section of the \KZ connection
\begin{equation}\label{eq:KZ}
\nabla\KKZ^{\hbbar}=
d-\hbbar\back\sum_{1\leq i<j\leq n}\back
\frac{dz_{i}-dz_{j}}{z_{i}-z_{j}}\Omega_{ij}\k
\end{equation}
iff the $V_{\lambda}\n[\mu]$--valued function
$\displaystyle{
g=f\cdot
\!\!\!\!\!\prod_{1\leq i<j\leq n}\!\!\!\!
(z_{i}-z_{j})^{h(\mu_{i}+\mu_{j}+2\mu_{i}\mu_{j}/k)}}$
is a horizontal section of
\begin{equation}\label{eq:Casimir}
\nabla_{\kappa}^{h}=
d-h\back\sum_{1\leq i<j\leq n}\back
\frac{dz_{i}-dz_{j}}{z_{i}-z_{j}}\kappa_{ij}\n
%\frac{\partial g}{\partial z_{i}}=
%h\sum_{j\neq i}\frac{\kappa_{ij}\n}{z_{i}-z_{j}}g
\end{equation}
where $\hbbar=2h$.
\end{theorem}
\proof Let $1\k=\sum_{i=1}^{k}E\k_{ii}$ be the generator
of the centre of $\gl_{k}$ so that, in obvious notation,
$\wt{\Omega}\k_{ij}=\Omega_{ij}\k+\frac{1}{k}\pi_{i}(1\k)
\pi_{j}(1\k)$. The operators $2\wt{\Omega}\k_{ij}$ and
$\kappa_{ij}\n$ both act on $M_{\lambda}^{\mu}\cong
V_{\lambda}\n[\mu]$ and, by proposition \ref{pr:Omegaij=Cij},
their restrictions differ by $-\mu_{i}-\mu_{j}$. The claim
follows since, for any $1\leq l\leq n$, $\pi_{l}(1\k)$
acts on $\Smu{}$ as multiplication by $\mu_{l}$ \halmos\\

\remark When $k=2$ and $\lambda$ is of the form $(|\mu|/2,
|\mu|/2,0,\lldots,0)$, where $|\mu|=\sum_{i=1}\mu_{i}$,
theorem \ref{th:KZ=KZ} is a representation--theoretic
analogue of the coincidence between the Kapovich--Millson
bending flows on the space of $n$--gons in $\IR^{3}$ with
side lengths $\mu_{1},\lldots,\mu_{n}$ \cite{KM} and the
Gel'fand--Cetlin flows on the Grassmannian Gr$_{2}(\IC^{n})$
\cite{GS} observed by Hausmann and Knutson in the context of
Gel'fand--McPherson duality \cite{HK}. I am grateful to
J. Millson for a careful explanation of this coincidence.\\

\remark An interesting relation between $\nablak$ and the
\KZ connection was recently noted by Felder, Markov, Tarasov
and Varchenko in \cite{FMTV}, where a variant of the connection
\eqref{eq:Casimir connection} is independently introduced and
studied. One of the main results of \cite{FMTV} is that, for
any simple Lie algebra $\g$, the connection $\nablak$ with
values in a tensor product $V_{1}\dots{\otimes}V_{n}$ of $n$
simple $\g$--modules is, when supplemented by suitable dynamical
parameters, bispectral to (\ie commutes with) the \KZ connection
for $\g$ with values in the same $n$--fold tensor product. An
analogous result is obtained in \cite{TV} for a difference
analogue of the connection $\nablak$. By comparison, theorem
\ref{th:KZ=KZ} can only hold for $\g=\sl_{n}$,
since it relies on the 'coincidence' of the regular Cartan of
$\gl_{n}$ with the configuration space of $n$ ordered points
in $\IC$, and asserts the {\it equality} of the two connections.\\

To relate the monodromy representations of $B_{n}$
corresponding to $\nabla\KKZ^{\hbbar}$ and $\nabla_{\kappa}
^{h}$, we need to specify how these induce flat connections
on $X_{n}/\SS_{n}$ and $X_{n}^{0}/\SS_{n}$ respectively.
For $\nabla\KKZ^{\hbbar}$, we let $\SS_{n}$ act on the
fibre
\begin{equation}
\bigoplus_{\nu\in\SS_{n}\mu}M_{\lambda}^{\nu}
\subset
\bigoplus_{\nu\in\SS_{n}\mu}{\mathcal S}^{\nu}\IC^{k}
\end{equation}
by permuting the tensor factors and take the quotient
connection. For $\nabla_{\kappa}^{h}$, we use the
construction of proposition \ref{th:existence}
and the homomorphism $\sigma:B_{n}\rightarrow
SL_{n}(\IC)$ given by
\begin{equation}\label{eq:assign}
T_{j}\longrightarrow
\exp(E_{j,j+1}\n)\exp(-E_{j+1,j}\n)\exp(E_{j,j+1}\n)=
\left(\begin{array}{ccrrcccc}
1&      & &  & & &      & \\  
 &\ddots& &  & & &      & \\
 &      &1&  & & &      & \\
 &      & &0 &1& &      & \\
 &      & &-1&0& &      & \\
 &      & &  & &1&      & \\
 &      & &  & & &\ddots& \\
 &      & &  & & &      &1
\end{array}
\right)
\end{equation}
where the off--diagonal terms are the $(j,j+1)$ and
$(j+1,j)$ entries. A direct computation, or \cite
[thm. 3.3]{Ti}, show that the assignement \eqref
{eq:assign} does indeed extend to a homomorphism
$B_{n}\rightarrow SL_{n}(\IC)$. Choose the generators
$T_{1}\ldots T_{n-1}$ of $B_{n}$ as at the end of
section \ref{se:flat}. 

\begin{corollary}\label{co:KZ Howe}
Let $\mu$ be a weight of $V_{\lambda}$ and
\begin{xalignat}{2}
\pi_{\kappa}^{h}&:
B_{n}\rightarrow
GL(\bigoplus_{\nu\in\SN{n}\mu}V_{\lambda}[\nu]),&
\pi\KKZ^{\hbbar}&:
B_{n}\rightarrow
GL(\bigoplus_{\nu\in\SN{n}\mu}M_{\lambda}^{\nu})
\end{xalignat}
the monodromy representations of the braid group $B_{n}$
corresponding to the connections \eqref{eq:Casimir} and
\eqref{eq:KZ} respectively. Then, for any $j=1\ldots n-1$,
\begin{equation}\label{eq:KZ=T}
\pi\KKZ^{2h}(T_{j})=\pi_{\kappa}^{h}(T_{j})
\cdot
e^{-\pi ih(E\n_{jj}+E\n_{j+1j+1}+2E\n_{jj}E\n_{j+1j+1}/k)}
\cdot
e^{i\pi E\n_{jj}}
\end{equation}
\end{corollary}
\proof
Let $s_{j}\in SL_{n}(\IC)$ be the right--hand side of
\eqref{eq:assign} so that $s_{j}=(j\medspace j+1)\cdot
e^{i\pi E_{jj}\n}$ in $GL_{n}(\IC)$. Let $\P\KKZ^{\hbbar}$,
$\P_{\kappa}^{h}$ denote parallel transport for $\nabla\KKZ
^{\hbbar}$ and $\nabla^{h}_{\kappa}$ respectively.
Then, by theorem \ref{th:KZ=KZ} and proposition
\ref{pr:recipe}, the following holds on
$M_{\lambda}^{\nu}\cong V_{\lambda}[\nu]$,
\begin{equation}
\begin{split}
\pi\KKZ^{2h}(T_{j})
&=
(j\medspace j+1)\P\KKZ^{2h}(T_{j})\\[1.2 ex]
&=
(j\medspace j+1)
e^{-\pi ih(\nu_{j}+\nu_{j+1}+2\nu_{j}\nu_{j+1}/k)}
\P_{\kappa}^{h}(T_{j})\\[1.2 ex]
&=
s_{j}\P_{\kappa}^{h}(T_{j})
e^{i\pi E_{jj}\n}
e^{-\pi ih(
E_{jj}\n+E_{j+1j+1}\n+2E_{jj}\n E_{j+1j+1}\n/k)}\\[1.2 ex]
&=
\pi_{\kappa}^{h}(T_{j})
e^{-\pi ih(
E_{jj}\n+E_{j+1j+1}\n+2E_{jj}\n E_{j+1j+1}\n/k)}
e^{i\pi E_{jj}\n}
\end{split}
\end{equation}
as claimed \halmos

\section{The quantum group $\Uhgl{p}$}\label{se:Uhg}
%=====================================

In this, and the following sections, we work over the
ring $\ICh$ of formal power series in the variable $\hbar$.
All tensor products of $\ICh$--modules are understood
to be completed in the $\hbar$--adic topology. For
$p\in\IN$, let $a_{ij}=2\delta_{ij}-\delta_{|i-j|=1}$,
$1\leq i,j\leq p$, be the entries of the Cartan matrix of
type $A_{p-1}$ and let $\Uhgl{p}$ be the corresponding
Drinfeld--Jimbo quantum group \cite{Dr1,Ji1} \ie the algebra
over $\ICh$ topologically generated by elements $E_{i},
F_{i}$, $i=1\ldots p-1$ and $D_{i}$, $i=1\ldots p$ subject
to the $q$--Serre relations
\begin{gather}
[D_{i},D_{j}]=0
%[H_{i},H_{j}]=0
\label{eq:Serre 1}\\[1.1 ex]
[D_{i},E_{j}]= (\delta_{ij}-\delta_{ij+1})E_{j}
%[H_{i},E_{j}]= a_{ij}E_{j}
\qquad
[D_{i},F_{j}]=-(\delta_{ij}-\delta_{ij+1})F_{j}
%[H_{i},F_{j}]=-a_{ij}F_{j}
\label{eq:Serre 2}\\[1.1 ex]
[E_{i},F_{j}]=\delta_{ij}
\frac{e^{\hbar H_{i}}-e^{-\hbar H_{i}}}
{e^{\hbar}-e^{-\hbar}}
%[E_{i},F_{j}]=\delta_{ij}
%\frac{e^{\hbar H_{i}}-e^{-\hbar H_{i}}}{e^{\hbar}-e^{-\hbar}}
\label{eq:Serre 3}\\[1.1 ex]
\sum_{k=0}^{1-a_{ij}}(-1)^{k}
\bin{1-a_{ij}}{k}
E_{i}^{k}E_{j}E_{i}^{1-a_{ij}-k}=0,
\quad\forall i\neq j
\label{eq:Serre 4}\\[1.1 ex]
\sum_{k=0}^{1-a_{ij}}(-1)^{k}
\bin{1-a_{ij}}{k}
F_{i}^{k}F_{j}F_{i}^{1-a_{ij}-k}=0,
\quad\forall i\neq j
\label{eq:Serre 5}
\end{gather}
where $H_{i}=D_{i}-D_{i+1}$ and, for any $n\geq k\in\IN$, 
\begin{align}
[n] &=\frac{e^{n\hbar}-e^{-n\hbar}}
          {e^{\hbar}-e^{\hbar}}\\
[n]!&=[n][n-1]\cdots[1]\\
\bin{n}{k} &= \frac{[n]!}{[k]![n-k]!}
\end{align}

$\Uhgl{p}$ is a topological Hopf algebra with coproduct
$\Delta$ and counit $\varepsilon$ given by
\begin{gather}
\Delta(D_{i})=D_{i}\otimes 1+1\otimes D_{i}
%\Delta(H_{i})=H_{i}\otimes 1+1\otimes H_{i}
\label{eq:coprod 1}\\
\Delta(E_{i})=E_{i}\otimes e^{\hbar H_{i}}+1\otimes E_{i}
\label{eq:coprod 2}\\
\Delta(F_{i})=F_{i}\otimes 1+e^{-\hbar H_{i}}\otimes F_{i}
\label{eq:coprod 3}
\end{gather}
and
\begin{equation}
\varepsilon(E_{i})=\varepsilon(F_{i})=\varepsilon(D_{i})=0
\end{equation}
Note that $I=D_{1}+\cdots+D_{p}$ is central so
that $\Uhgl{p}\cong\Uhsl{p}\otimes\IC[I]\fml$
as Hopf algebras where the coproduct on $\IC[I]
\fml$ is given by $\Delta(I)=I\otimes 1+1\otimes
I$ and $\Uhsl{p}\subset\Uhgl{p}$ is the closed
Hopf subalgebra generated by $E_{i},F_{i}$ and
$H_{i}$, $i=1\ldots p-1$.\\

By a finite--dimensional representation of
$\Uhgl{p}$ we shall mean a $\Uhgl{p}$--module
which is topologically free and finitely
generated over $\ICh$ and on which $I$ acts
semisimply with eigenvalues in $\IC$.
Choose an algebra isomorphism $\phi:\Uhgl{p}
\rightarrow\Ugl{p}\fml$ mapping each $D_{i}$
onto $E_{ii}$ \cite[prop. 4.3]{Dr2} and let
$V$ be a finite--dimensional $\gl_{p}$--module
on which $1\p=\sum_{i=1}^{p}E_{ii}$ acts semisimply.
Then, $\Ugl{p}\fml$ acts on $V\fml$ and the
latter becomes, via $\phi$, a finite--dimensional
representation of $\Uhgl{p}$. Conversely,

\begin{proposition}\label{pr:at 0}
Let $\V$ be a finite--dimensional representation
of $\Uhgl{p}$ and $V=\V/\hbar V$ the corresponding
$\gl_{p}$--module. Then, as $\Uhgl{p}$--modules,
\begin{equation}
\V\cong V\fml
\end{equation}
\end{proposition}
\proof
Since $I$ is diagonalisable on $\V$ and commutes
with $\Uhgl{p}$, we may assume that it acts on
$\V$ as multiplication by a scalar $\lambda\in\IC$.
Since $\V$ is topologically free, $\V\cong V\fml$
as $\ICh$--modules so that $\V$ is a deformation
of the finite--dimensional $\sl_{p}$--module $V$.
Since $\sl_{p}$ is simple, $H^{1}(\sl_{p},V)=0$
and $\V$ is isomorphic, as $\sl_{p}$, and therefore
as $\gl_{p}$--module to the trivial deformation
of $V$. Thus, $\V\cong V\fml$ as $U\gl_{p}\fml$,
and therefore as $\Uhgl{p}$--modules \halmos

\begin{corollary}\label{co:mult h}
Let $U,V$ be finite--dimensional $\gl_{p}$--modules
on which $1\p$ acts semisimply. If $U\otimes V$
decomposes as
\begin{equation}\label{eq:mult}
U\otimes V\cong\bigoplus_{W}N_{W}W
\end{equation}
for some $\gl_{p}$--modules $W$ and multiplicities
$N_{W}\in\IN$, then, as $\Uhgl{p}$--modules,
\begin{equation}\label{eq:mult h}
U\fml\otimes V\fml\cong\bigoplus_{W}N_{W}W\fml
\end{equation}
\end{corollary}
\proof
By \eqref{eq:mult}, both sides of \eqref{eq:mult h}
have the same specialisation at $\hbar=0$ and are
therefore isomorphic by proposition \ref{pr:at 0}
\halmos

\section{The dual pair $(\Uhgl{k},\Uhgl{n})$}\label{se:q Howe}
%============================================

We shall need the analogue of theorem \ref{th:Howe}
in the setting of the algebra $\SMh$ of functions
on quantum $k\times n$ matrix space. With the exception
of theorems \ref{th:pre q Howe} and \ref{th:qHowe},
this section follows \cite[\S 1.5]{Ba} (see also
\cite{Ga}). By definition, $\Sh{k}{n}$ is the
algebra over $\ICh$ topologically generated by
elements $X_{ij}$, $1\leq i\leq k$, $1\leq j\leq
n$ with relations 
\begin{equation}\label{eq:Manin}
X_{ij}X_{kl}=\left\{\begin{array}{cl}
X_{kl}X_{ij}&
\text{if $k>i$ and $l<j$ or $k<i$ and $l>j$}\\[1.4ex]
e^{-\hbar}X_{kl}X_{ij}&
\text{if $k>i$ and $l=j$ or $k=i$ and $l>j$}\\[1.4ex]
X_{kl}X_{ij}-(e^{\hbar}-e^{-\hbar})X_{kj}X_{il}
&\text{if $k>i$ and $l>j$}
\end{array}\right.
\end{equation}
$\Sh{k}{n}$ is $\IN$--graded by decreeing that
each $X_{ij}$ is of degree 1 and we denote its
homogeneous components by $\SMhd$, $d\in\IN$.
For any $k\times n$ matrix $m$ with entries
$m_{ij}\in\IN$, set
\begin{equation}
\begin{split}
\XX^{m}
&=
X_{11}^{m_{11}}\cdots X_{k1}^{m_{k1}}
\cdots
X_{1n}^{m_{1n}}\cdots X_{kn}^{m_{kn}}\\
&=
X_{11}^{m_{11}}\cdots X_{1n}^{m_{1n}}
\cdots
X_{k1}^{m_{k1}}\cdots X_{kn}^{m_{kn}}
\end{split}
\end{equation}
%be the corresponding column and row monomials.
By the commutation relations \eqref{eq:Manin},
the $\XX^{m}$, with $m\in\M_{k,n}$ such that
$|m|=d$, span $\SMhd$, where $|m|=\sum_{i,j}
m_{ij}$.

\begin{theorem}[Parshall--Wang]\label{th:basis}
The monomials $\XX^{m}$, $m\in{\M}_{kn}(\IN)$,
are linearly independent over $\ICh$. In particular,
the set $\{\XX^{m}\}_{|m|=d}$ is a $\ICh$--basis
of $\SMhd$.
\end{theorem}
\proof
This is proved in \cite[thm. 3.5.1]{PW} for $k=n$
and over the field $\IC(q)$ of rational functions
of $q=e^{\hbar}$ rather than over $\ICh$.
The proof however works equally well for $k\neq n$
and, as remarked in \cite{PW}, over $\ICh$
\halmos\\

As in the classical case, $\SMh$ is a module algebra
over $\Uhgl{k}\otimes\Uhgl{n}$. This may be seen in
the following way. For any $l\in\IN$, one readily
checks that the assignement
\begin{equation}
X_{ij}\rightarrow\sum_{l'=1}^{l}X_{il'}\otimes X_{l'j}
\end{equation}
extends uniquely to an algebra homomorphism
$\Delta_{kln}:\Sh{k}{n}\rightarrow\Sh{k}{l}\otimes\Sh{l}{n}$
such that, for any $l,m\in\IN$, the following
diagram commutes
\begin{equation}\label{eq:Delta klm}
\begin{diagram}
\Sh{k}{n}                &\rTo^{\Delta_{kln}}         &\Sh{k}{l}\otimes\Sh{l}{n}\\
\dTo^{\Delta_{kmn}}      &                            &\dTo_{1\otimes\Delta_{lmn}}\\
\Sh{k}{m}\otimes\Sh{m}{n}&\rTo_{\Delta_{klm}\otimes 1}&\Sh{k}{l}\otimes\Sh{l}{m}\otimes\Sh{m}{n}\\
\end{diagram}
\end{equation}
In particular, $\Sh{k}{k}$ and $\Sh{n}{n}$ are topological
bialgebras with comultiplications $\Delta_{kkk}$
and $\Delta_{nnn}$ respectively and counit $\varepsilon
(X_{ij})=\delta_{ij}$. Moreover, the maps $\Delta_{kkn}$
and $\Delta_{knn}$ give $\Sh{k}{n}$ the structure of
a $\Sh{k}{k}$--$\Sh{n}{n}$ bicomodule algebra each
homogeneous component of which is invariant under
$\Sh{k}{k}$ and $\Sh{n}{n}$ since $\Delta_{kln}
(\SMhdp{k}{n})\subset\SMhdp{k}{l}\otimes\SMhdp{l}{n}$.\\

We shall need a columnwise (resp. rowwise) description
of the coaction of $\Sh{k}{k}$ (resp. $\Sh{n}{n}$) on
$\Sh{k}{n}$. Consider the quantum $k$ and $n$--dimensional
planes \ie the algebras $\Sh{k}{1}$ and $\Sh{1}{n}$.
By the commutation relations \eqref{eq:Manin} and
theorem \ref{th:basis}, these may be embedded as
subalgebras of $\Sh{k}{n}$ via the maps
\begin{gather}
c_{j}:\Sh{k}{1}\rightarrow\Sh{k}{n},
\quad
c_{j}(X_{i1})=X_{ij}\\
r_{i}:\Sh{1}{n}\rightarrow\Sh{k}{n},
\quad
r_{i}(X_{1j})=X_{ij}
\end{gather}
with $1\leq i\leq k,1\leq j\leq n$. By \eqref{eq:Delta klm},
$\Sh{k}{1}$ is a left algebra comodule over $\Sh{k}{k}$ and
$\Sh{1}{n}$ a right algebra comodule over $\Sh{n}{n}$.

\begin{lemma}\label{le:qtensor}
As left, $\IN$--graded $\Sh{k}{k}$--comodules,
\begin{equation}\label{eq:Phi}
\Sh{k}{n}\cong\Sh{k}{1}^{\otimes n}
\end{equation}
via the map $\Phi:p_{1}\dots{\otimes}p_{n}
\rightarrow c_{1}(p_{1})\cdots c_{n}(p_{n})$.
Similarly, as right, $\IN$--graded
$\Sh{n}{n}$--comodules,
\begin{equation}\label{eq:Psi}
\Sh{k}{n}\cong\Sh{1}{n}^{\otimes k}
\end{equation}
via $\Psi:q_{1}\dots{\otimes}q_{k}\rightarrow
r_{1}(q_{1})\cdots r_{k}(q_{k})$.
\end{lemma}
\proof
The map $\Phi$ clearly preserves the grading and,
by theorem \ref{th:basis}, restricts to a $\ICh
$--linear isomorphism of homogeneous components
since it bijectively maps the monomial basis of
$\Sh{k}{1}^{\otimes n}$ onto the basis $\XX^{m}$
of $\Sh{k}{n}$. The fact that $\Phi$ itself is an
isomorphism follows easily because any element of
$\Sh{1}{n}^{\otimes k}$ or $\Sh{k}{n}$ is the
convergent sum of its homogeneous components.
As readily checked,
the diagram
\begin{equation}
\begin{diagram}
\Sh{k}{k}\otimes\Sh{k}{1}&\rTo^{1\otimes c_{j}}&\Sh{k}{k}\otimes\Sh{k}{n}\\
\uTo^{\Delta_{kk1}}      &                     &\uTo_{\Delta_{kkn}}\\
\Sh{k}{1}                &\rTo_{c_{j}}         &\Sh{k}{n}
\end{diagram}
\end{equation}
is commutative for any $1\leq j\leq n$ and therefore
so is
\begin{equation}
\begin{diagram}
(\Sh{k}{k}\otimes\Sh{k}{1})^{\otimes n}
                           &\rTo^{\mu^{(n)}\otimes\Phi}&\Sh{k}{k}\otimes\Sh{k}{n}\\
\uTo^{\Delta_{kk1}^{\otimes n}}
                           &                           &\uTo_{\Delta_{kkn}}\\
\Sh{k}{1}^{\otimes n}      &\rTo_{\Phi}                &\Sh{k}{n}
\end{diagram}
\end{equation}
where $\mu^{(n)}:\Sh{k}{k}^{\otimes n}\rightarrow\Sh{k}
{k}$ is the $n$--fold multiplication. This proves
\eqref{eq:Phi}. The proof of \eqref{eq:Psi} is
identical \halmos\\

We turn now to the action of $\Uhgl{k}\otimes\Uhgl{n}$
on $\Sh{k}{n}$. For $p=k,n$, consider the vector
representation of $\Uhgl{p}$ \ie the module $V=\IC^{p}
\fml$ with basis $e_{1},\ldots e_{p}$ and action
given by
\begin{equation}\label{eq:vector}
D_{i}=E_{ii},\qquad
E_{i}=E_{ii+1},\qquad
F_{i}=E_{i+1i}
\end{equation}
where $E_{ab}e_{c}=\delta_{bc}e_{a}$.
Let $e^{1},\ldots,e^{p}\in V^{*}$ be the dual basis
of $e_{1},\ldots,e_{p}$, $\Uhgl{p}\res$ the restricted
dual of $\Uhgl{p}$  and $t_{ij}\in\Uhgl{p}\res$
the matrix coefficient defined by
\begin{equation}
t_{ij}(x)=\<e^{i},xe_{j}\>
\end{equation}

\begin{proposition}\label{pr:pairing}
The assignement $X_{ij}\rightarrow t_{ij}$ extends
uniquely to a bialgebra morphism $\kappa_{p}:\Sh{p}
{p}\rightarrow\Uhgl{p}\res$.
\end{proposition}
\proof
We need to check that the $t_{ij}$ satisfy
the relations \eqref{eq:Manin}, \ie that
when evaluated on $\Delta(x)$, $x\in\Uhgl{p}$,
\begin{equation}\label{eq:t Manin}
t_{ij}\otimes t_{kl}=\left\{\begin{array}{cl}
t_{kl}\otimes t_{ij}&
\text{if $k>i$ and $l<j$ or $k<i$ and $l>j$}\\[1.4ex]
e^{-\hbar}t_{kl}\otimes t_{ij}&
\text{if $k>i$ and $l=j$ or $k=i$ and $l>j$}\\[1.4ex]
t_{kl}\otimes t_{ij}-(e^{\hbar}-e^{-\hbar})t_{kj}\otimes t_{il}
&\text{if $k>i$ and $l>j$}
\end{array}\right.
\end{equation}
Let $R^{\vee}=\sigma\cdot R\in\End(V\otimes V)$
where $\sigma\in GL(V\otimes V)$ is the
flip and $R$ is the universal $R$--matrix
of $\Uhgl{p}$ acting on $V\otimes V$. Then
\cite{Ji2},\cite[\S 8.3.G]{CP}
\begin{equation}
R=\left(
%e^{-\hbar/p}
e^{\hbar}\sum_{i=1}^{p}E_{ii}\otimes E_{ii}+
\sum_{1\leq i\neq j\leq p}E_{ii}\otimes E_{jj}+
(e^{\hbar}-e^{-\hbar})\sum_{1\leq i<j\leq p}
E_{ij}\otimes E_{ji}\right)
\end{equation}
so that the matrix entries of $R^{\vee}$ are
\begin{equation}\label{eq:R entries}
R^{\vee}_{ik,jl}=
%e^{-\hbar/p}\cdot
\left\{\begin{array}{cl}
e^{\delta_{jl}\hbar}
&\text{if $i=l$ and $k=j$}\\[1.4ex]
e^{\hbar}-e^{-\hbar}
&\text{if $i=j$, $k=l$ and $j>l$}\\[1.4ex]
0&\text{otherwise}
\end{array}\right.
\end{equation}
From \eqref{eq:R entries}, one readily checks
that both sides of \eqref{eq:t Manin} coincide
when evaluated on any $A\in\End(V\otimes V)$
commuting with $R^{\vee}$. Since $R^{\vee}$ is
a $\Uhgl{p}$--intertwiner, $\kappa$ extends to
an algebra morphism which respects the counit
and coproduct since $\Delta(t_{ij})=\sum_{q=1}
^{p}t_{iq}\otimes t_{qj}$ \halmos

\begin{theorem}\label{th:pre q Howe}\hfill
\begin{enumerate}
\item The maps $\kappa_{p}$, $p=k,n$ of proposition
\ref{pr:pairing} give $\Sh{k}{n}$ the structure of
an algebra module over $\Uhgl{k}\otimes\Uhgl{n}$ with
invariant homogeneous components $\SMhd$, $d\in\IN$.
\item The maps $\Phi,\Psi$ of lemma \ref{le:qtensor}
yield isomorphisms
\begin{equation}
\Sh{k}{n}\cong\Sh{k}{1}^{\otimes n}
\qquad\text{and}\qquad
\Sh{k}{n}\cong\Sh{1}{n}^{\otimes k}
\end{equation}
as $\IN$--graded $\Uhgl{k}$ and $\Uhgl{n}$--modules
respectively.
\item The action of the generators $E\p_{q},F\p_{q}$,
$q=1\ldots p-1$ and $D\p_{q}$, $q=1\ldots p$ of $\Uhgl
{p}$, $p=k,n$ in the monomial basis $\XX^{m}$, $m\in
\MM{k}{n}$, is given by
\begin{align}
D_{i}\k\medspace\XX^{m}&=
\sum_{j=1}^{n}m_{ij}\medspace\XX^{m}
%H_{i}\k\medspace\XX^{m}&=
%\sum_{j=1}^{n}(m_{ij}-m_{i+1j})\medspace\XX^{m}
\label{eq:explicit first}\\
E_{i}\k\medspace\XX^{m}&=
\sum_{j=1}^{n}[m_{i+1j}]
\prod_{j'=j+1}^{n}e^{\hbar(m_{ij'}-m_{i+1j'})}
\medspace\XX^{m+\varepsilon_{ij}-\varepsilon_{i+1j}}
\label{eq:explicit second}\\
F_{i}\k\medspace\XX^{m}&=
\sum_{j=1}^{n}[m_{ij}]
\prod_{j'=1}^{j-1}e^{-\hbar(m_{ij'}-m_{i+1j'})}
\medspace\XX^{m-\varepsilon_{ij}+\varepsilon_{i+1j}}
\label{eq:explicit middle}
\end{align}
%for $1\leq i\leq k-1$,
where $(\varepsilon_{ab})_{cd}=\delta_{ac}\delta_{bd}$,
and
\begin{align}
D_{j}\n\medspace\XX^{m}&=
\sum_{i=1}^{k}m_{ij}\medspace\XX^{m}
%H_{j}\n\medspace\XX^{m}&=
%\sum_{i=1}^{k}(m_{ij}-m_{ij+1})\medspace\XX^{m}
\label{eq:explicit middle'}\\
E_{j}\n\medspace\XX^{m}&=
\sum_{i=1}^{k}[m_{ij+1}]
\prod_{i'=i+1}^{k}e^{\hbar(m_{i'j}-m_{i'j+1})}
\medspace\XX^{m+\varepsilon_{ij}-\varepsilon_{ij+1}}
\label{eq:explicit second'}\\
F_{j}\n\medspace\XX^{m}&=
\sum_{i=1}^{k}[m_{ij}]
\prod_{i'=1}^{i-1}e^{-\hbar(m_{i'j}-m_{i'j+1})}
\medspace\XX^{m-\varepsilon_{ij}+\varepsilon_{ij+1}}
\label{eq:explicit last}
\end{align}
%for $1\leq j\leq n-1$.
\end{enumerate}
\end{theorem}
\proof
Using the transposition anti--involution $\tau$
on $\Sh{k}{k}$ given by $\tau(X_{ij})=X_{ji}$, we may
regard $\Sh{k}{n}$ as a right algebra comodule over
$\Sh{k}{k}\otimes\Sh{n}{n}$ and therefore, via the
pairings $\<\cdot,\cdot\>:\Sh{m}{m}\otimes\Uhgl{m}
\rightarrow\ICh$, $m=k,n$, given by proposition
\ref{pr:pairing} as a left algebra module over
$\Uhgl{k}\otimes\Uhgl{n}$. This proves (i) and (ii).
Explicitly, for $x^{(m)}\in\Uhgl{m}$, $m=k,n$ and
$p\in\Sh{k}{n}$
\begin{align}
x\k\medspace p&=
\<x\k\otimes 1,\tau\otimes 1\cdot\Delta_{kkn}(p)\>
\label{eq:left}\\
x\n\medspace p&=\<1\otimes x\n,\Delta_{knn}(p)\>
\label{eq:right}
\end{align}
Using \eqref{eq:left} and \eqref{eq:vector}, one gets
\begin{align}
D\k_{i}\medspace X_{i'j}&=
\delta_{ii'}\medspace X_{ij}\\
%H\k_{i}\medspace X_{i'j}&=
%(\delta_{ii'}-\delta_{i+1i'})\medspace X_{i'j}\\
E\k_{i}\medspace X_{i'j}&=
\delta_{i+1 i'}\medspace X_{ij}\\
F\k_{i}\medspace X_{i'j}&=
\delta_{ii'}\medspace X_{i+1j}
\end{align}
Using the algebra module property $x(pq)=\mu(\Delta(x)
p\otimes q)$ where $x\in\Uhgl{k}$, $p,q\in\Sh{k}{n}$
and $\mu:\Sh{k}{n}^{\otimes 2}\rightarrow\Sh{k}{n}$ is
multiplication, \eqref{eq:coprod 1}--\eqref{eq:coprod 3}
and the commutation relations \eqref{eq:Manin} shows
by induction on $m\in\IN$ that
\begin{align}
D\k_{i}\medspace X_{i'j}^{m}&=
\delta_{ii'}m\medspace X_{ij}^{m}\\
%H\k_{i}\medspace X_{i'j}^{m}&=
%(\delta_{ii'}-\delta_{i+1i'})m\medspace X_{i'j}^{m}\\
E\k_{i}\medspace X_{i'j}^{m}&=
\delta_{i+1 i'}[m]\medspace X_{ij}X_{i+1j}^{m-1}\\
F\k_{i}\medspace X_{i'j}^{m}&=
\delta_{ii'}[m]\medspace X_{ij}^{m-1}X_{i+1j}
\end{align}
Let $\Delta^{(a)}:\Uhgl{k}\rightarrow\Uhgl{k}^
{\otimes a}$, $a\in\IN^{*}$ be recursively defined
by $\Delta^{(1)}=\id$, $\Delta^{(a+1)}=\Delta\otimes
\id^{\otimes(a-1)}\cdot\Delta^{(a)}$. Then, by
\eqref{eq:coprod 1}--\eqref{eq:coprod 3}
\begin{align}
\Delta^{(a)}D\k_{i}&=\sum_{b=1}^{a}
1^{\otimes(b-1)}\otimes D\k_{i}\otimes 1^{\otimes(a-b)}
%\Delta^{(a)}H\k_{i}&=\sum_{b=1}^{a}
%1^{\otimes(b-1)}\otimes H\k_{i}\otimes 1^{\otimes(a-b)}
\label{eq:a coprod 1}\\
\Delta^{(a)}E\k_{i}&=\sum_{b=1}^{a}
1^{\otimes(b-1)}\otimes E\k_{i}\otimes(e^{\hbar H\k_{i}})^{\otimes(a-b)}\\
\Delta^{(a)}F\k_{i}&=\sum_{b=1}^{a}
(e^{-\hbar H\k_{i}})^{\otimes(b-1)}\otimes F\k_{i}\otimes 1^{\otimes(a-b)}
\label{eq:a coprod 3}
\end{align}
The formulae \eqref{eq:explicit first}--\eqref{eq:explicit middle}
now follow from the algebra module property and
\eqref{eq:a coprod 1}--\eqref{eq:a coprod 3}.
The proof of \eqref{eq:explicit middle'}--\eqref{eq:explicit last}
is similar \halmos\\

The following result is proved in \cite{Ba} and \cite{Ga}
for the quantum groups $U_{q}\gl_{k},U_{q}\gl_{n}$ by a
different method

\begin{theorem}\label{th:qHowe}
For any $d\in\IN$, the $\Uhgl{k}\otimes\Uhgl{n}$--module
$\SMhd$ decomposes as
\begin{equation}
\SMhd\cong
\bigoplus_{\substack{\lambda\in\IY_{\min(k,n)}\\|\lambda|=d}}
V\k_{\lambda}\fml{\otimes}V\n_{\lambda}\fml
\end{equation}
\end{theorem}
\proof
By theorem \ref{th:basis}, $\SMhd$ has no torsion, and
is therefore a topologically free $\ICh$--module.
Moreover, by \eqref{eq:explicit first}--\eqref{eq:explicit last},
$\SMhd/\hbar\SMhd$ is the $\gl_{k}\oplus\gl_{n}$--module
$\SMd$. The conclusion follows from theorem \ref{th:Howe}
and proposition \ref{pr:at 0} \halmos

\section{Braid group actions on quantum matrix space}
\label{se:Bn on Mkn}
%====================================================

We compare in this section two actions of the
braid group $B_{n}$ on the algebra $\Sh{k}{n}$
of functions of quantum $k\times n$ matrix
space. The first is the $R$--matrix representation
obtained by regarding $\Sh{k}{n}$ as the $\Uhgl
{k}$--module $\Sh{k}{1}^{\otimes n}$. The second
is the quantum Weyl group action of $B_{n}$ on
$\Sh{k}{n}$ viewed as a $\Uhgl{n}$--module. We
will show that these representations essentially
coincide, thus extending to the $q$--setting the
fact that the symmetric group $\mathfrak{S}_{n}$
acts on $(S^{\bullet}\IC^{k})^{\otimes n}\cong\SM$
via the permutation matrices in $GL_{n}(\IC)$.\\

More precisely, for any $1\leq j\leq n$, let
$R_{j}^{\vee}$ be the universal $R$--matrix
of $\Uhgl{k}$ acting on the $j$ and $j+1$
tensor copies of $\Sh{k}{n}\cong\Sh{k}{1}^
{\otimes n}$ and $S_{j}$ the quantum Weyl
group element of $\Uhgl{n}$ corresponding
to the simple root $\alpha_{j}=\theta_{j}-
\theta_{j+1}$. We will show that

\begin{equation}\label{eq:Sj Rj}
R_{j}^{\vee}=S_{j}
\cdot e^{-\hbar(D_{j}\n+D_{j}\n D_{j+1}\n/k)}
\cdot e^{i\pi D_{j}\n}
\end{equation}

where $D_{1}\n,\ldots,D_{n}\n$ are the generators
of the Cartan subalgebra of $\Uhgl{n}$. The
proof of \eqref{eq:Sj Rj} is based upon the
following observation, which we owe to B. Feigin.
Both sides of \eqref{eq:Sj Rj} only act upon the
$j$ and $j+1$ tensor copies of $\Sh{k}{1}^
{\otimes n}$ so that its proof reduces to
a computation in $\Sh{k}{1}^{\otimes 2}
\cong\Sh{k}{2}$. Since both sides intertwine
the action of $\Uhgl{k}$ on $\Sh{k}{2}$, it
suffices to compare them on highest weight vectors.
These, and the action of $R_{j}^{\vee}$ are computed
in \S \ref{ss:R on singular}. The action of
$S_{j}$ is computed in \S \ref{ss:q Weyl on singular}.\\

\remark It is easy to check that neither action
of $B_{n}$ is compatible with the algebra structure
of $\Sh{k} {n}$, so that \eqref{eq:Sj Rj} cannot
be proved by merely checking it on the generators
$X_{ij}$. This stems from the fact that quantum
Weyl group operators are not group--like.

\subsection{$R$--matrix action on singular vectors}
\label{ss:R on singular}

For any $d\in\IN$, let $\Shd$ be the homogeneous
component of degree $d$ of $\Sh{k}{1}$. By \eqref
{eq:explicit first}--\eqref{eq:explicit middle},
$\Shd$ is a deformation of the $d$--th symmetric
power $\Sd$ of the vector representation of $\gl_{k}$.
Let $\mu_{1},\mu_{2}\in\IN$, then

\begin{lemma}\label{le:highest weight}
As $\Uh\gl_{k}$--modules,
\begin{equation}\label{eq:q-Pieri}
\Shmu{1}\otimes\Shmu{2}\cong
\bigoplus_{i=0}^{\min(\mu_{1},\mu_{2})}
V_{\young{\mu_{1}+\mu_{2}-i}{i}}\k\fml
\end{equation}
where $V_{\young{a}{b}}\k$ is the irreducible
representation of $\gl_{k}$ with highest weight
$(a,b,0,\ldots,0)$. The corresponding highest weight
vectors $v^{\mu_{1},\mu_{2}}_{i}$ are given by
%and $v\op_{\mu_{1},\mu_{2}}(i)$ for $\Delta$
%and $\Delta\op$ respectively 
\begin{equation}\label{eq:v(i)}
v^{\mu_{1},\mu_{2}}_{i}=
\sum_{a=0}^{i}(-1)^{a}\bin{i}{a}e^{\hbar a(\mu_{2}-a+1)}
X_{11}^{\mu_{1}-i+a}X_{21}^{i-a}\otimes
X_{12}^{\mu_{2}-a}X_{22}^{a}
\end{equation}
%v\op_{\mu_{1},\mu_{2}}(i)&=
%\sum_{a=0}^{i}(-1)^{a}\bin{i}{a}e^{\hbar a(\mu_{1}-a+1)}
%X_{11}^{\mu_{1}-a}X_{21}^{a}\otimes
%X_{12}^{\mu_{2}-i+a}X_{22}^{i-a}
%\label{eq:v(i)opp}
\end{lemma}
\proof The decomposition \eqref{eq:q-Pieri} follows
from the Pieri rules for $\gl_{k}$ and corollary
\ref{co:mult h}. Fix $i\in\{0,\ldots,\min(\mu_{1},
\mu_{2})\}$. By \eqref{eq:explicit first}, any
$v\in\Shmu{1}\otimes\Shmu{2}$ of weight $(\mu_{1}
+\mu_{2}-i,i,0\ldots,0)$ is of the form
\begin{equation}
v=\sum_{a=0}^{i}c_{a}\medspace
X_{11}^{\mu_{1}-i+a}X_{21}^{i-a}\otimes
X_{12}^{\mu_{2}-a}X_{22}^{a}
\end{equation}
for some constants $c_{a}\in\IC$. By
\eqref{eq:coprod 2} and \eqref{eq:explicit second},
$\Delta(E_{j})v=0$ for any $j\geq 2$ so that $v$
is a highest weight vector iff
\begin{equation}
\begin{split}
\Delta(E_{1})v
&=
\sum_{a=0}^{i}
c_{a}[i-a]e^{\hbar(\mu_{2}-2a)}
X_{11}^{\mu_{1}-i+a+1}X_{21}^{i-a-1}\otimes
X_{12}^{\mu_{2}-a}X_{22}^{a}\\
&+
\sum_{a=0}^{i}
c_{a}[a]
X_{11}^{\mu_{1}-i+a}X_{21}^{i-a}\otimes
X_{12}^{\mu_{2}-a+1}X_{22}^{a-1}
\end{split}
\end{equation}
is equal to zero. This yields $c_{a}=-c_{a-1}\medspace
e^{\hbar(\mu_{2}-2a+2)}[i-a+1]/[a]$ and therefore
\begin{equation}
c_{a}=(-1)^{a}\bin{i}{a}e^{\hbar a(\mu_{2}-a+1)}c_{0}
\end{equation}
whence \eqref{eq:v(i)} \halmos\\
%\eqref{eq:v(i)opp} follows from \eqref{eq:v(i)} by
%permuting $\mu_{1}$ and $\mu_{2}$

Let $R$ be the universal $R$--matrix of $\Uhsl{k}$
and $R^\vee=\sigma\cdot R:\Shmu{1}\otimes\Shmu{2}
\rightarrow\Shmu{2}\otimes\Shmu{1}$ the corresponding
$\Uhsl{k}$--intertwiner, where $\sigma$ is the
permutation of the tensor factors.

\begin{proposition}\label{pr:explicit Rv}
The following holds on $\Shmu{1}\otimes\Shmu{2}
\bigoplus\Shmu{2}\otimes\Shmu{1}$,
\begin{equation}\label{eq:Rv(i)}
R^{\vee}\medspace v^{\mu_{1},\mu_{2}}_{i}=
(-1)^{i}
e^{\hbar((\mu_{1}-i)(\mu_{2}-i)-i-\mu_{1}\mu_{2}/k)}
\medspace
v^{\mu_{2},\mu_{1}}_{i}
\end{equation}
%R\medspace v\op_{\mu_{1},\mu_{2}}(i)&=
%(-1)^{i}e^{\hbar\left(
%(\mu_{1}-i)(\mu_{2}-i)-i-\frac{\mu_{1}\mu_{2}}{k}
%\right)}
%v_{\mu_{1},\mu_{2}}(i)
%\label{eq:Rv(i)opp}
\end{proposition}
\proof
For any $1\leq i\leq k-1$, let $s_{i}=(i\medspace
i+1)\in\SS_{k}$ be the $i$th elementary transposition
and let
\begin{equation}
w_{0}=
(1\medspace k)(2\medspace k-1)
\cdots
(\lfloor\frac{k}{2}\rfloor\medspace\lceil\frac{k}{2}\rceil)
\end{equation}
be the longest element of $\SS_{k}$. Consider the
following reduced expression for $w_{0}$
\begin{equation}
w_{0}=
s_{k-1}\cdots s_{1}
s_{k-1}\cdots s_{2}
\cdots\cdots
s_{k-1}s_{k-2}s_{k-1}
=s_{i_{1}}\cdots s_{i_{k(k-1)/2}}
\end{equation}
and let $\beta_{j}=s_{i_{1}}\cdots s_{i_{j-1}}
(\theta_{i_{j}}-\theta_{i_{j}+1})$ be the associated
enumeration of the positive roots of $\sl_{k}$ so that
\newcommand {\spazio}{\negthickspace\negthickspace\negthinspace}
\begin{equation}
\begin{array}{lllrl}
 \beta_{1}=\theta_{k-1}-\theta_{k},
&\beta_{2}=\theta_{k-2}-\theta_{k},
&\cdots
&\beta_{k-1}&\spazio=\theta_{1}-\theta_{k},\\
&\beta_{k}=\theta_{k-2}-\theta_{k-1},
&\cdots
&\beta_{2k-3}&\spazio=\theta_{1}-\theta_{k-1},\\
&&\ddots&&\\%\vdots\\
&&&\beta_{k(k-1)/2}&\spazio=\theta_{1}-\theta_{2}
\end{array}
\end{equation}
Let $E_{\beta_{j}},F_{\beta_{j}}\in\Uhsl{k}$, $1\leq
j\leq k(k-1)/2$ be the corresponding quantum root vectors
so that $E_{\beta_{j}}=E_{i}$ and $F_{\beta_{j}}=F_{i}$
whenever $\beta_{j}$ is the simple root $\alpha_{i}$
\cite[prop. 1.8]{Lu2}. Then, \cite{KR,LS,Ro},\cite[thm.
8.3.9]{CP}
\begin{equation}\label{eq:R matrix}
R=\exp\left(\hbar\sum_{i=1}^{k-1}H^{i}\otimes H_{i}\right)
\prod_{j=1}^{k(k-1)/2}
\exp_{q}\left((q-q^{-1})E_{\beta_{j}}\otimes F_{\beta_{j}}\right)
\end{equation}
where $\{H^{i}\}_{i=1}^{k-1}\subset\h$ is the basis
of the Cartan subalgebra of $\sl_{k}$ dual to $\{H_{i}
\}_{i=1}^{k-1}$ with respect to the pairing $\<X,Y\>
=\tr(XY)$, $q=e^{\hbar}$,
\begin{equation}\label{eq:qexp}
\exp_{q}(x)=\sum_{n\geq 0} q^{n(n-1)/2}\frac{x^{n}}{[n]!}
\end{equation}
and the product in \eqref{eq:R matrix} is taken so
that the factor $\exp_{q}\left((q-q^{-1})E_{\beta_{j}}
\otimes F_{\beta_{j}}\right)$ is placed to the left
of $\exp_{q}\left((q-q^{-1})E_{\beta_{j'}}\otimes
F_{\beta_{j'}}\right)$ whenever $j>j'$. To compute
$R\medspace v^{\mu_{1},\mu_{2}}_{i}$, note that for
any positive root $\beta\neq\theta_{1}-\theta_{2}$
and $0\leq a\leq\mu_{1}$,
\begin{equation}
E_{\beta}\medspace X_{11}^{\mu_{1}-a}X_{21}^{a}=0
\end{equation}
since, by \eqref{eq:explicit first}, $(\mu_{1}-a,a,
0,\ldots,0)+\beta$ is not a weight of $\Shmu{1}$.
Thus, using
\eqref{eq:explicit first}--\eqref{eq:explicit middle}
\begin{equation}
\begin{split}
R\medspace v^{\mu_{1},\mu_{2}}_{i}
&=
\exp\left(\hbar\sum_{i=1}^{k-1}H^{i}\otimes H_{i}\right)
\exp_{q}\left((q-q^{-1})E_{1}\otimes F_{1}\right)
v^{\mu_{1},\mu_{2}}_{i}\\
&=
\sum_{\substack{0\leq a\leq i\\[1.1ex]0\leq n\leq i-a}}
(-1)^{a}e^{\hbar a(\mu_{2}-a+1)}\bin{i}{a}
e^{\hbar((\mu_{1}-i+a+n)(\mu_{2}-a-n)+(i-a-n)(a+n)-\mu_{1}\mu_{2}/k)}\\[1.2ex]
&\cdot
\frac{e^{\hbar n(n-1)/2}(e^{\hbar}-e^{-\hbar})^n}{[n]!}
\frac{[i-a]![\mu_{2}-a]!}{[i-a-n]![\mu_{2}-a-n]!}\\[1.2ex]
&\cdot
X_{11}^{\mu_{1}-i+a+n}X_{21}^{i-a-n}\otimes
X_{12}^{\mu_{2}-a-n}X_{22}^{a+n}
\end{split}
\end{equation}
which, upon setting $\alpha=a+n$, yields
\begin{equation}
e^{\hbar((\mu_{1}-i)(\mu_{2}-i)-i-\mu_{1}\mu_{2}/k)}
\sum_{\alpha=0}^{i}(-1)^{\alpha}\bin{i}{\alpha}
e^{\hbar(i-\alpha)(\mu_{1}-i+\alpha+1)}
S^{\mu_{2}}_{\alpha}\medspace
X_{11}^{\mu_{1}-i+\alpha}X_{21}^{i-\alpha}
\otimes
X_{12}^{\mu_{2}-\alpha}X_{22}^{\alpha}
\end{equation}
where
\begin{equation}
S^{\mu}_{\alpha}=
e^{\hbar\alpha(\mu-\alpha+1)}
\sum_{n=0}^{\alpha}(-1)^{n}\bin{\alpha}{n}
\frac{[\mu-\alpha+n]!}{[\mu-\alpha]!}
e^{\hbar(\alpha-n)(\mu-\alpha+n+1)+\hbar n(n-1)/2}
(e^{\hbar}-e^{-\hbar})^{n}
\end{equation}
We claim that $S^{\mu}_{\alpha}=1$ for any $\alpha
\leq\mu\in\IN$ so that \eqref{eq:Rv(i)} holds.
Indeed, using
\begin{equation}
\bin{\alpha}{a}=
\delta_{\alpha>a}\bin{\alpha-1}{a}e^{-\hbar a}+
\delta_{a>0}\bin{\alpha-1}{a-1}e^{\hbar(\alpha-a)}
\end{equation}
one readily finds
\begin{equation}
S_{\alpha}^{\mu}=
e^{2\hbar(\mu-\alpha+1)}S_{\alpha-1}^{\mu-1}-
(e^{2\hbar(\mu-\alpha+1)}-1)S^{\mu}_{\alpha-1}
\end{equation}
%\begin{equation}
%\begin{split}
%S_{\alpha}^{\mu}
%&=
%e^{2\hbar\alpha(\mu-\alpha+1)}
%\sum_{a=0}^{\alpha-1}(-1)^{a}\bin{\alpha-1}{a}
%\frac{[\mu-\alpha+a]!}{[\mu-\alpha]!}
%(1-e^{-2\hbar})^{a}
%e^{-\hbar a}e^{-\hbar a(2\mu-4\alpha+a+1)/2}\\[1.1 ex]
%&-
%e^{2\hbar\alpha(\mu-\alpha+1)}
%\sum_{a=0}^{\alpha-1}(-1)^{a}\bin{\alpha-1}{a}
%\frac{[\mu-\alpha+a+1]!}{[\mu-\alpha]!}
%(1-e^{-2\hbar})^{a+1}
%e^{\hbar(\alpha-a-1)}e^{-\hbar(a+1)
%(2\mu-4\alpha+a+2)/2}\\[1.1 ex]
%&=
%e^{2\hbar(\mu-\alpha+1)}S_{\alpha-1}^{\mu-1}-
%e^{\hbar(\mu-\alpha+2)}(1-e^{-2\hbar})
%[\mu-\alpha+1]S^{\mu}_{\alpha-1}\\[1.1 ex]
%&=
%e^{2\hbar(\mu-\alpha+1)}S_{\alpha-1}^{\mu-1}-
%(e^{2\hbar(\mu-\alpha+1)}-1)S^{\mu}_{\alpha-1}
%\end{split}
%\end{equation}
whence $S^{\mu}_{\alpha}=1$ by induction on $\alpha$
since $S^{\mu}_{0}=1$ for any $\mu\in\IN$ \halmos\\
%Finally, \eqref{eq:Rv(i)opp} follows from \eqref
%{eq:Rv(i)} by permuting $\mu_{1}$ and $\mu_{2}$ 

\subsection{Quantum Weyl group action on singular vectors}
\label{ss:q Weyl on singular}

Let $E,F,H$ be the standard generators of $\Uh\sl_{2}$.

\begin{lemma}\label{le:sl2 on singular}
The following holds in $\Sh{k}{2}$,
\begin{align}
E \medspace v^{\mu_{1},\mu_{2}}_{i}&=
[\mu_{2}-i] \medspace v^{\mu_{1}+1,\mu_{2}-1}_{i}\\
F \medspace v^{\mu_{1},\mu_{2}}_{i}&=
[\mu_{1}-i] \medspace v^{\mu_{1}-1,\mu_{2}+1}_{i}\\
H \medspace v^{\mu_{1},\mu_{2}}_{i}&=
(\mu_{1}-\mu_{2}) \medspace v^{\mu_{1},\mu_{2}}_{i}
\label{eq:H}
\end{align}
\end{lemma}
\proof By \eqref{eq:explicit second'},
\begin{equation}
\begin{split}
E\medspace v^{\mu_{1},\mu_{2}}_{i}
&=
\sum_{a=0}^{i}(-1)^{a}\bin{i}{a}e^{\hbar a(\mu_{2}-a+1)}
[\mu_{2}-a]e^{\hbar(i-2a)}
X_{11}^{\mu_{1}-i+a+1}X_{21}^{i-a}\otimes
X_{12}^{\mu_{2}-a-1}X_{22}^{a}\\
&+
\sum_{a=0}^{i}(-1)^{a}\bin{i}{a}e^{\hbar a(\mu_{2}-a+1)}
[a]
X_{11}^{\mu_{1}-i+a}X_{21}^{i-a+1}\otimes
X_{12}^{\mu_{2}-a}X_{22}^{a-1}\\
&=
\sum_{a=0}^{i}(-1)^{a}\bin{i}{a}
e^{\hbar a(\mu_{2}-a+1)}
[\mu_{2}-a]e^{\hbar(i-2a)}
X_{11}^{\mu_{1}-i+a+1}X_{21}^{i-a}\otimes
X_{12}^{\mu_{2}-a-1}X_{22}^{a}\\
&-
\sum_{a=0}^{i}(-1)^{a}\bin{i}{a}
e^{\hbar (a+1)(\mu_{2}-a)}[i-a]
X_{11}^{\mu_{1}-i+a+1}X_{21}^{i-a}\otimes
X_{12}^{\mu_{2}-a-1}X_{22}^{a}\\
&=
[\mu_{2}-i]
\sum_{a=0}^{i}(-1)^{a}\bin{i}{a}
e^{\hbar a(\mu_{2}-a)}
X_{11}^{\mu_{1}+1-i+a}X_{21}^{i-a}\otimes
X_{12}^{\mu_{2}-1-a}X_{22}^{a}\\
&=
[\mu_{2}-i]
v^{\mu_{1}+1,\mu_{2}-1}_{i}
\end{split}
\end{equation}
Similarly, by \eqref{eq:explicit last},
\begin{equation}
\begin{split}
F\medspace v^{\mu_{1},\mu_{2}}_{i}
&=
\sum_{a=0}^{i}(-1)^{a}\bin{i}{a}e^{\hbar a(\mu_{2}-a+1)}
[\mu_{1}-i+a]
X_{11}^{\mu_{1}-i+a-1}X_{21}^{i-a}\otimes
X_{12}^{\mu_{2}-a+1}X_{22}^{a}\\
&+
\sum_{a=0}^{i}(-1)^{a}\bin{i}{a}e^{\hbar a(\mu_{2}-a+1)}
[i-a]e^{-\hbar(\mu_{1}-\mu_{2}-i+2a)}
X_{11}^{\mu_{1}-i+a}X_{21}^{i-a-1}\otimes
X_{12}^{\mu_{2}-a}X_{22}^{a+1}\\
&=
\sum_{a=0}^{i}(-1)^{a}\bin{i}{a}e^{\hbar a(\mu_{2}-a+1)}
[\mu_{1}-i+a]
X_{11}^{\mu_{1}-i+a-1}X_{21}^{i-a}\otimes
X_{12}^{\mu_{2}-a+1}X_{22}^{a}\\
&-
\sum_{a=0}^{i}(-1)^{a}\bin{i}{a}
e^{\hbar (a-1)(\mu_{2}-a+2)}[a]
e^{-\hbar(\mu_{1}-\mu_{2}-i+2a-2)}
X_{11}^{\mu_{1}-i+a-1}X_{21}^{i-a}\otimes
X_{12}^{\mu_{2}-a+1}X_{22}^{a}\\
&=
[\mu_{1}-i]
\sum_{a=0}^{i}(-1)^{a}\bin{i}{a}
e^{\hbar a(\mu_{2}-a+2)}
X_{11}^{\mu_{1}-1-i+a}X_{21}^{i-a}\otimes
X_{12}^{\mu_{2}+1-a}X_{22}^{a}\\
&=
[\mu_{1}-i]
v^{\mu_{1}-1,\mu_{2}+1}_{i}
\end{split}
\end{equation}
Finally, \eqref{eq:H} follows from \eqref{eq:explicit middle'}
\halmos\\

Let now
\begin{equation}\label{eq:generator}
S=
\exp_{q^{-1}}(q^{-1}Eq^{-H})
\exp_{q^{-1}}(-F)
\exp_{q^{-1}}(qEq^{H})
q^{H(H+1)/2}
\end{equation}
be the generator of the quantum Weyl group of $\Uh\sl_{2}$
(\cite{Lu3,KR,So}, we use the form given in \cite{Ka,Sa})
where $q=e^{\hbar}$ and the $q$--exponential is defined by
\eqref{eq:qexp}.

\begin{proposition}\label{pr:S on singular}
The following holds in $\Sh{k}{2}$,
\begin{equation}\label{eq:S on singular}
S v^{\mu_{1},\mu_{2}}_{i}=
(-1)^{\mu_{1}-i}q^{(\mu_{1}-i)(\mu_{2}-i+1)}
v^{\mu_{2},\mu_{1}}_{i}
\end{equation}
\end{proposition}
\proof
Fix $\mu,i\in\IN$ with $2i\leq\mu$. By lemma
\ref{le:sl2 on singular}, the vectors $u_{k}
=v^{\mu-i-k,i+k}_{i}$, with $0\leq k\leq\mu-
2i$ satisfy
\begin{align}
E\medspace u_{k}&=[k]\medspace u_{k-1}\\
F\medspace u_{k}&=[\mu-2i-k]\medspace u_{k+1}\\
H\medspace u_{k}&=(\mu-2i-2k)\medspace u_{k}
\end{align}
and therefore span the indecomposable $\Uhsl
{2}$--module of dimension $\mu-2i+1$. Moreover,
\begin{equation}\label{eq:up down}
u_{k}=
\frac{[\mu-2i-k]!}{[\mu-2i]!}F^{k}u_{0}=
\frac{[\mu-2i-k]!}{[\mu-2i]!}E^{\mu-2i-k}u_{\mu-2i}
\end{equation}
Since $\Ad(S)H=-H$, $Su_{0}$ is proportional
to $u_{\mu-2i}$ and, by \eqref{eq:generator}
is therefore equal to
\begin{equation}
(-1)^{\mu-2i}q^{\mu-2i}
\frac{F^{\mu-2i}}{[\mu-2i]!}u_{0}
=
(-1)^{\mu}q^{\mu-2i}u_{\mu-2i}
\end{equation}
Next, using $\Ad(S)F=-q^{-H}E$ and \eqref{eq:up down},
we find
\begin{equation}
S\medspace u_{k}
=
\frac{[\mu-2i-k]!}{[\mu-2i]!}\Ad(S)F^{k}\medspace Su_{0}
=
(-1)^{\mu-k}q^{(k+1)(\mu-2i-k)}u_{\mu-2i-k}
\end{equation}
Thus, setting $\mu=\mu_{1}+\mu_{2}$, so that
$v^{\mu_{1},\mu_{2}}_{i}=u_{\mu_{2}-i}$, we find
\begin{equation}
S\medspace v^{\mu_{1},\mu_{2}}_{i}
=
(-1)^{\mu_{1}-i}
q^{(\mu_{1}-i)(\mu_{2}-i+1)}v^{\mu_{2},\mu_{1}}_{i}
\end{equation}
as claimed \halmos

\subsection{Identification of R and quantum Weyl group actions}
\label{ss:R=S}

Fix $1\leq j\leq n$ and let $R_{j}^\vee$ be the universal
$R$--matrix of $\Uh\gl_{k}$ acting on the $j$ and $j+1$
tensor copies of $\Sh{k}{1}^{\otimes n}$. Let $S_{j}$
be the element of the quantum Weyl group of $\Uh\gl_{n}$
corresponding to the simple root $\theta_{j}-\theta_{j+1}$.

\begin{theorem}\label{th:R=S}
The following holds on $\Sh{k}{n}$
\begin{equation}\label{eq:R=S}
R_{j}^{\vee}=S_{j}
\cdot e^{-\hbar(D_{j}\n+D_{j}\n D_{j+1}\n/k)}
\cdot e^{i\pi D_{j}\n}
\end{equation}
\end{theorem}
\proof
Let $\Uhgl{2}^{j}\subset\Uhgl{n}$ be the Hopf subalgebra
generated by $E_{j}\n,F_{j}\n,D_{j}\n,D_{j+1}\n$. By
\eqref{eq:explicit middle'}--\eqref{eq:explicit last},
$\Uhgl{2}^{j}$ only acts upon the variables $X_{ij},
X_{ij+1}$, $1\leq i\leq k$. Thus, $\Uhgl{2}^{j}$ and
therefore $S_{j}$ only act on the $j$ and $j+1$ tensor
copies $\Sh{k}{1}_{j}$, $\Sh{k}{1}_{j+1}$ of $\Sh{k}{1}
^{\otimes n}\cong\Sh{k}{n}$. Since this is also the
case of $R_{j}^{\vee}$, the proof of \eqref{eq:R=S}
reduces to a computation on the $\Uhgl{k}\otimes\Uhgl{2}
^{j}$--module $\Sh{k}{1}_{j}\otimes\Sh{k}{1}_{j+1}
\cong\Sh{k}{2}$. Both sides of \eqref{eq:R=S}
clearly commute with $\Uhgl{k}$, so it is suffices to
check their equality on the singular vectors of
\begin{equation}
\Sh{k}{1}_{j}\otimes\Sh{k}{1}_{j+1}=
\bigoplus_{\mu_{1},\mu_{2}\in\IN}
\Shmu{1}\otimes\Shmu{2}
\end{equation}
\ie on the vectors $v^{\mu_{1},\mu_{2}}_{i}\in\Shmu{1}
\otimes\Shmu{2}$ of lemma \ref{le:highest weight}.
By propositions \ref{pr:explicit Rv} and
\ref{pr:S on singular},
\begin{equation}
R_{j}^{\vee}\medspace v^{\mu_{1},\mu_{2}}_{i}
=
S_{j}\cdot e^{-\hbar(\mu_{1}+\mu_{1}\mu_{2}/k)}(-1)^{\mu_{1}}
\medspace v^{\mu_{1},\mu_{2}}_{i}
\end{equation}
whence \eqref{eq:R=S} since, for any $1\leq l\leq
n$, $D_{l}\n$ gives the $\IN$--grading on $\Sh{k}
{1}_{l}$ \halmos\\

\remark The coincidence of the two representations of
$B_{n}$ studied above was also noted by Baumann \cite
[prop. 12]{Ba} in the special case when both actions
are restricted to the subspace $\Shmu{1}\dots{\otimes}
\Shmu{n}$ of $\SMh$ where $\mu_{1}=\cdots=\mu_{n}=1$.

\section{Monodromic realisation of quantum Weyl group operators}
\label{se:main}
%===============================================================

The following is the main result of this paper. It
was conjectured for any simple Lie algebra $\g$ by
De Concini around 1995 (unpublished) and, independently,
in \cite{TL}. We prove it here for $\g=\sl_{n}(\IC)$.

\begin{theorem}\label{th:main}
Let $\g=\sl_{n}(\IC)$, $V$ a finite--dimensional $\g$--module,
$\mu$ a weight of $V$ and
\begin{equation}
V^{\mu}=\bigoplus_{\nu\in W\mu}V[\nu]
\end{equation}
the direct sum of the weight spaces of $V$ corresponding
to the Weyl group orbit of $\mu$. Let $\pi_{\kappa}^{h}:
\Bg\rightarrow GL(V^{\mu}\fmll)$ be the monodromy representation
of the connection
\begin{equation}
d-h
\sum_{\alpha>0}\frac{d\alpha}{\alpha}
\medspace\pi_{V}(\kappa_{\alpha})
%d-h\back\sum_{1\leq i<j\leq n}\back
%\frac{d(z_{i}-z_{j})}{z_{i}-z_{j}}\pi_{V}(\kappa_{ij})
\end{equation}
obtained by regarding $h$ as a formal variable. Let
$\pi_{W}:\Bg\rightarrow GL(V^{\mu}\fmll)$ be the quantum
Weyl group action obtained by regarding $V\fmll$ as a
$\Uhg$--module with $\hbar=2\pi ih$. Then, $\pi_{\kappa}^{h}$
and $\pi_{W}$ are equivalent.
\end{theorem}
\proof
We may assume that $V$ is irreducible with highest weight
$\lambda=(\lambda_{1},\ldots,\lambda_{n})$ where $\lambda
_{i}\in\IN$. Regard $V$ as a $\gl_{n}$--module by letting
$1\n=\sum_{j=1}^{n}E_{jj}\n$ act as multiplication by $|
\lambda|=\sum_{i=1}^{n}\lambda_{i}$ and fix some $k\geq n$.
By lemma \ref{le:singular=weight}, $V[\nu]$ is isomorphic
to the space $M_{\lambda}^{\nu}$ of singular vectors of
weight $\lambda$ for the diagonal action of $\gl_{k}$ on
\begin{equation}
\SmC{\nu}=
\SmC{\nu_{1}}\dots{\otimes}\SmC{\nu_{n}}\subset
\SM
\end{equation}
and, by corollary \ref{co:KZ Howe}, the monodromy
representation of the \KZ connection \eqref{eq:KZ} on
$\bigoplus_{\nu\in\SS_{k}}M_{\lambda}^{\nu}$ and $\pi_
{\kappa}^{h}$ are related by
\begin{equation}
\pi\KKZ^{\hbbar}(T_{j})=\pi_{\kappa}^{h}(T_{j})
e^{-\pi ih(E\n_{jj}+E\n_{j+1 j+1}+2E\n_{jj}E\n_{j+1 j+1}/k)}
e^{i\pi E\n_{jj}}
\end{equation}
where $\hbbar=2h$.\\

We shall now use the Kohno--Drinfeld theorem to relate
$\pi^{\hbbar}\KKZ$ to the $R$--matrix representation
of $B_{n}$ corresponding to the action of $\Uhgl{k}$
on $\SmC{\nu}\fml$. Let for this purpose $\phi:\Uhgl
{k}\rightarrow U\gl_{k}\fml$ be an algebra isomorphism
whose reduction mod $\hbar$ is the identity and which
acts as the identity on the Cartan subalgebras \ie $
\phi(D_{i}\k)=E_{ii}\k$, $1\leq i\leq k$ \cite[prop.
4.3]{Dr2}. Then, $\Uhgl{k}$ act on each $\SmC{\nu_{j}}
\fml$ via $\phi$ and on
\begin{equation}
\SmC{\nu_{1}}\fml\dots{\otimes}\SmC{\nu_{n}}\fml=
\SmC{\nu_{1}}\dots{\otimes}\SmC{\nu_{n}}\fml=
\SmC{\nu}\fml
\end{equation}
via the $n$--fold coproduct $\Delta\n:\Uhgl{k}\rightarrow
\Uhgl{k}^{\otimes n}$ recursively defined by
\begin{align}
\Delta^{(1)}   &=
\id
\label{eq:nfold 1}\\
\Delta^{(a+1)} &=
\Delta\otimes\id^{\otimes(a-1)}\circ\Delta^{(a)},
\quad a\geq 1
\label{eq:nfold 2}
\end{align}
where $\Delta=\Delta^{(2)}$ is given by
\eqref{eq:coprod 1}--\eqref{eq:coprod 3}.\\

Let $\Delta_{0}$ be the standard, cocommutative coproduct on
$U\gl_{k}$ so that $U\gl_{k}$ acts on $\SmC{\nu}$ via $\Delta_{0}\n:U\gl
_{k}\rightarrow U\gl_{k}^{\otimes(n)}$ defined as in \eqref
{eq:nfold 1}--\eqref{eq:nfold 2} with $\Delta$ replaced by
$\Delta_{0}$. Since $\Delta=\Delta_{0}+o(\hbar)$ and $H^{1}
(\sl_{k},\sl_{k}\oplus\sl_{k})=0$, there exists a twist $F
=1\otimes 1+o(\hbar)\in U\gl_{k}^{\otimes(2)}\fml$ such
that, for any $x\in\Uhgl{k}$,
\begin{equation}
\phi\otimes\phi\circ\Delta(x)=F\Delta_{0}(\phi(x))F^{-1}
\end{equation}
It follows that the actions of $\Uhgl{k}$ and $U\gl_{k}\fml$
on $\SmC{\nu}$ are related by
\begin{equation}\label{eq:ad F(n)}
\phi^{\otimes(n)}\circ\Delta^{(n)}(x)=
F\n\Delta_{0}^{(n)}(\phi(x)){F\n}^{-1},
\quad x\in\Uhgl{k}
\end{equation}
where $F\n\in U\gl_{k}^{\otimes n}\fml$ is recursively
defined by
\begin{align}
F^{(1)}  &=
1\\
F^{(a+1)}&=
F\otimes 1^{\otimes(a-1)}\cdot
\Delta_{0}\otimes\id^{\otimes(a-1)}(F^{(a)}),
\quad a\geq 1
\end{align}

Let now $M_{\lambda}^{\nu}\subset\SmC{\nu}$ and
$M_{\hbar,\lambda}^{\nu}\subset\SmC{\nu}\fml$ be the
subspaces of vectors of highest weight $\lambda$ for
the actions of $\gl_{k}$ and $\Uhgl{k}$ respectively.
We shall need the following

\begin{lemma}
$\displaystyle
{F\n\medspace M_{\lambda}^{\nu}\fml=M_{\hbar,\lambda}^{\nu}}$.
\end{lemma}
\proof
Let $\SmC{\nu}(\lambda)\subset\SmC{\nu}$ and
$\SmC{\nu}\fml(\lambda)\subset\SmC{\nu}\fml$
be the isotypical components of types $V_{\lambda}$
and $V_{\lambda}\fml$ for $\gl_{k}$ and $\Uhgl{k}$
respectively. The subspace $F\n\SmC{\nu}(\lambda)\fml$
is invariant under $\Uhgl{k}$ by \eqref{eq:ad F(n)}
and its reduction mod $\hbar$ is equal to $\SmC
{\nu}(\lambda)$ since $F\n=1^{\otimes(n)}+o(\hbar)$.
Thus, by proposition \ref{pr:at 0}, $F\n\SmC{\nu}
(\lambda)\fml$ is isomorphic to $\SmC{\nu}(\lambda)
\fml$ and therefore contained in $\SmC{\nu}\fml(
\lambda)$. Since this holds for any $\lambda$,
the inclusion is an equality. Noting now that
$F\n$ is equivariant for the action of the
Cartan subalgebras of $\gl_{k}$ and $\Uhgl{k}$
since $\Delta_{0}\circ\phi(h)=\phi\otimes\phi\circ
\Delta(h)$ for any $h\in\h$, we get that
\begin{equation}
\begin{split}
F\n\medspace M_{\lambda}^{\nu}\fml
&=
F\n\medspace\left(\SmC{\nu}(\lambda)[\lambda]\fml\right)\\
&=
F\n\medspace\left(\SmC{\nu}(\lambda)\fml\right)[\lambda]\\
&=
\SmC{\nu}\fml(\lambda)[\lambda]\\
&=
M_{\hbar,\lambda}^{\nu}
\end{split}
\end{equation}
as claimed \halmos\\

Summarising, we have an action of $B_{n}$ on $\bigoplus
_{\nu\in\SS_{n}\mu}M_{\lambda}^{\nu}$ via the monodromy
of the \KZ connection and an action of $B_{n}$ on
$\bigoplus_{\nu\in\SS_{n}\mu}M_{\hbar,\lambda}^{\nu}$
via the $R$--matrix representation of $\Uhsl{k}$.
Drinfeld's theorem \cite{Dr3,Dr4,Dr5} asserts that
the twist $F$ may be chosen so that
\begin{equation}
F\n:
\bigoplus_{\nu\in\SS_{n}\mu}M_{\lambda}^{\nu}\fml
\longrightarrow
\bigoplus_{\nu\in\SS_{n}\mu}
M_{\hbar,\lambda}^{\nu}
\end{equation}
is $B_{n}$--equivariant so that, for any $1\leq j\leq n-1$,
\begin{equation}
\pi\KKZ^{\hbbar}(T_{j})={F\n}^{-1}R_{j}^{\vee}{F\n}
\end{equation}
where $R_{j}^{\vee}$ is the universal $R$--matrix
for $\Uhsl{k}$ acting on the $j$ and $j+1$ copies
of $\bigoplus_{\nu\in\SS_{n}\mu}\SmC{\nu}\fml$ and
$\hbar=\pi i\hbbar$.\\

Let now $\Sh{k}{n}\cong(\Sh{k}{1})^{\otimes n}$ be the
algebra of functions on quantum $k\times n$ matrix space
defined in section \ref{se:q Howe} and, for any $\nu\in
\IN^{n}$, let
\begin{equation}
\SmCh{\nu}=\SmCh{\nu_{1}}\dots{\otimes}\SmCh{\nu_{n}}
\end{equation}
be the space of polynomials which are homogeneous of
degree $\nu_{j}$ in the variables $X_{1j},\ldots,X_{kj}$,
for any $1\leq j\leq n$. By
\eqref{eq:explicit first}--\eqref{eq:explicit middle},
$\SmCh{\nu}/\hbar\SmCh{\nu}$ is the $\gl_{k}$--module
$\SmC{\nu}$ so that, by proposition \ref{pr:at 0}, we
may identify $\SmCh{\nu}$ with $\SmC{\nu}\fml$ as
$\Uhgl{k}$--modules.\\

By theorem \ref{th:qHowe} and the fact that the
$\IN$--grading on the $j$th tensor factor of $\Sh
{k}{n}\cong(\Sh{1}{k})^{\otimes n}$ is given by
the action of $D_{j}\n$, the space $M_{\hbar,\lambda}^{\nu}$
is isomorphic to the subspace $V_{\lambda}\n[\nu]
\fml\subset V_{\lambda}\n\fml$ of weight $\nu$.
Using now theorem \ref{th:R=S}, we find
\begin{equation}\label{eq:almost}
\begin{split}
\pi_{\kappa}^{h}(T_{j})
&=
{F\n}^{-1} \cdot S_{j}
\cdot e^{-\hbar(D_{j}\n+D_{j}\n D_{j+1}\n/k)}
\cdot e^{i\pi D_{j}\n}\cdot {F\n}\\
&\cdot
e^{\pi ih(E\n_{j}+E\n_{j+1}+2E\n_{j}E\n_{j+1}/k)}
\cdot e^{i\pi E\n_{j}}
\end{split}
\end{equation}
where $S_{j}=\pi_{W}(T_{j})$ is the quantum Weyl
group element of $\Uhgl{n}$ corresponding to the
simple root $\alpha_{j}=\theta_{j}-\theta_{j+1}$.
Since for any $l$, $E_{ll}\n=D_{l}\n$ on $V^{\mu}
\cong\bigoplus_{\nu\in\SS_{k}}M_{\lambda}^{\nu}$
we get
\begin{equation}
\begin{split}
\pi_{\kappa}^{h}(T_{j})
&=
{F\n}^{-1} S_{j} e^{-\hbar/2(D_{j}\n-D_{j+1}\n)} {F\n}\\
&=
{F\n}^{-1} S_{j} e^{-\hbar H_{j}\n/2} {F\n}\\
&=
{F\n}^{-1} e^{\hbar\rho\n/2}S_{j} e^{-\hbar\rho\n/2} {F\n}
\end{split}
\end{equation}
where $\rho\n=\half{1}\sum_{j=1}(n-2j+1)D_{j}\n$
is the half--sum of the positive (co)roots of
$\sl_{n}$ \halmos\\

\remark Theorem \ref{th:main} is proved in \cite{TL}
for the following pairs $(\g,V)$ where $\g$ is
a simple Lie algebra and $V$ an irreducible,
finite--dimensional representation
\begin{enumerate}
\item $\g=\sl_{2}$ and $V$ is any irreducible representation.
\item $\g=\sl_{n}$ and $V$ is a fundamental representation.
\item $\g=\so_{n}$ and $V$ is the vector or a spin representation.
\item $\g=\sp_{n}$ and $V$ is the defining vector representation.
\item $\g=\e_{6},\e_{7}$ and $V$ is a minuscule representation.
\item $\g={\mathfrak g}_{2}$ and $V$ is the 7--dimensional representation.
\item $\g$ is any simple Lie algebra and $V\cong\g$ its adjoint representation.
\end{enumerate}

\end{document}